\documentclass[10pt]{article}
\usepackage{hyperref}%,showlabels}
\usepackage{fancyvrb}  % verbatim
\usepackage{amsmath,amssymb,color,bbm,ifthen,enumerate}
\usepackage{graphicx}
\usepackage{multirow}

\newtheorem{thm}{Theorem}[section]

\newtheorem{prop}[thm]{Proposition}

\newenvironment{rem}[1][Remark.]{\begin{trivlist}
\item[\hskip \labelsep {\bfseries #1}]}{\end{trivlist}}
\newenvironment{rems}[1][Remarks.]{\begin{trivlist}
\item[\hskip \labelsep {\bfseries #1}]}{\end{trivlist}}

\newenvironment{pf}{\paragraph{Proof.}}
{\nopagebreak\hfill\nopagebreak\rule{2mm}{2mm}\par\bigskip}

\renewcommand{\hat}{\widehat}

\newcommand{\bone}{\mbox{\boldmath$1$}}

\begin{document}

\renewcommand{\baselinestretch}{1.3}

\title{Numerical Computation of First-Passage Times of Increasing L\'{e}vy
Processes}

\author{ Mark Veillette  and Murad S. Taqqu  \thanks{ This research was partially supported by the NSF grants DMS-0505747,
DMS-0706786, and DGE-0221680.}
 \\
Boston University 
}

%\address{ Dept. of Mathematics \& Statistics 
%111 Cummington St. 
% Boston, MA, 02215 } 

\maketitle

\begin{abstract}

Let $\{D(s), \ s \geq 0\}$ be a non-decreasing L\'{e}vy process.  The
first-hitting time process $\{E(t) \ t \geq 0\}$ (which is sometimes
referred to as an inverse subordinator) defined by $E(t) =
\inf \{s: D(s) > t \}$ is a process which has arisen in many
applications.  Of particular interest is the
mean first-hitting time $U(t)=\mathbb{E}E(t)$.  This function characterizes all
finite-dimensional distributions of the process $E$.  The function $U$
can be calculated by inverting the Laplace transform of the function
$\widetilde{U}(\lambda) = (\lambda \phi(\lambda))^{-1}$, where $\phi$
is the  L\'{e}vy exponent of the subordinator $D$. In this paper,
we give two methods for computing numerically the inverse of this Laplace
transform.  The first is based on the Bromwich integral
and the second is based on the Post-Widder inversion formula.  The
software written to support this work is available from the authors and we
illustrate its use at the end of the paper.

\end{abstract}

\section{Introduction}

Let $ \{D(s),\ s \geq 0 \} $ be a L\'{e}vy subordinator, that is, a
non-decreasing L\'{e}vy process
starting from $0$, which is continuous from the right with left
limits.  This process has stationary
and independent increments and is characterized by its Laplace Transform
\begin{equation}
 \mathbb{E} e^{-\lambda D(s)} = e^{-s \phi(\lambda)}, \quad \lambda
 \geq 0.
\end{equation}
The function $\phi$ above is known as the L\'{e}vy exponent (or
Laplace exponent) and is given by the L\'{e}vy-Khintchine formula:
\begin{equation}\label{e:LKformula}
\phi(\lambda) = \mu \lambda + \int_{(0,\infty)} \left( 1 - e^{-\lambda x} \right) \Pi(dx),
\end{equation}
where $\mu \geq 0$ is the drift and the L\'{e}vy measure $\Pi$ is a measure on
$\mathbb{R}^+ \cup \{0 \}$ which satisfies $\int_0^\infty (1 \wedge x)
\Pi(dx) < \infty$ (see ~\cite{Applebaum:2004},
~\cite{Bertoin:1996} or ~\cite{Bertoin:1999}). 

Consider the \textit{inverse subordinator} $\{E(t), \ t \geq 0 \}$, which is given by the
first-passage time of $D$:
\begin{equation}
E(t) = \inf\{s: D(s) > t \}.
\end{equation}
The process $E(t)$ is non-decreasing, and its sample paths are a.s. continuous if and only
if $D$ is strictly increasing.  Also, $E$ is, in general, non-Markovian with
non-stationary and non-independent increments.

Inverse subordinators are of particular interest in the study of
fractional kinetics and the scaling limits of continuous time
random walks, \cite{Baule:2005},~\cite{Baule:2007},~\cite{Meerschaert:2006},~\cite{Meerschaert:2004},
  ~\cite{Arous:2007}.   Here, a Markov process, $\{X(t), t \geq 0\}$, is time-changed with an
  inverse subordinator, yielding the new process $M(t) = X(E(t))$, $t \geq
  0$.   For certain choices of subordinators $D$, using $E$ as a time
  change gives a model of
  \textit{anomalous diffusion}, where the mean-squared displacement of
  $M$ grows non-linearly in time. For example, if $D$ is an
  $\alpha$-stable subordinator, that is, the
  subordinator whose L\'{e}vy exponent is given by
\begin{equation}\label{e:astable}
\phi(\lambda) = \lambda^\alpha, \quad 0 < \alpha < 1,
\end{equation}
then the mean of $E(t)$ is given by the power law
$\mathbb{E}E(t) \sim t^\alpha$.  Other types of non-linear behavior are also
possible.  In \cite{Meerschaert:2006}, Meerschaert considers the
subordinator with L\'{e}vy exponent given by the following
generalization of (\ref{e:astable}):
\begin{equation}\label{e:mixedintro}
\phi(\lambda) = \int_0^1 \lambda^\beta dp(\beta) ,
\end{equation}
where $p$ is a probability measure on $(0,1)$.  Depending on the
choice of $p$, the mean of $E$ in this case can grow at logarithmic
rates.  Meerschaert also shows in this case that the transition density of the
process $M$ solves a time-fractional diffusion equation whose
order is distributed according to $p$.

Inverse Subordinators have also appeared in many other areas of probability
theory.  Early work regarding the joint distribution of $E(t)$ and
$D(E(t))$ was done in \cite{Gusak:1969} and \cite{Kesten:1969}.  More
recently, in \cite{Kaj:2005}, Kaj and Martin-L\"{o}f show that a scaled sum of these processes converges
weakly to another non-stable and non-Gaussian process.  An application of
inverse subordinators in modeling of foreign exchange markets was
considered in \cite{Winkel:2005}.   Distributional properties of
inverse subordinators were studied in \cite{Nord:2005}, which drew
upon a connection with Cox processes. A different approach using differential equations to
characterize the joint distribution function of the process $E$ was
used in  \cite{veillette:2008} and
\cite{Baule:2005}.

 An important function in the study of inverse subordinators is the
 so-called \textit{renewal function}, $U(t)$, $t \geq 0$, which is given by the mean of the inverse subordinator, $U(t) =
\mathbb{E}E(t)$.  It has be shown by different methods that this
function characterizes all finite-dimensional distributions of the
process $E$.  The Laplace transform of $U$ is given simply in terms of
the L\'{e}vy exponent $\phi$, however, inverting this Laplace transform
to obtain $U$ in closed form is not always possible.  In this paper,
we provide two numerical methods for obtaining this mean first-passage
time which take as input the drift $\mu$ and the L\'{e}vy measure
$\Pi$ of the corresponding subordinator $D$.  Once $U$ is obtained,
higher order moments can be calculated by methods described in
\cite{veillette:2008} or \cite{Nord:2005}.

This paper is organized as follows: We begin in Section \ref{s:background} by giving a brief
background on the theory of inverse subordinators, as well as describe
the importance of the function $U$.  In Section \ref{s:numerical}, we
develop two methods for calculating $U$ numerically and test them in cases where $U$ can be
computed analytically.  In Section \ref{s:examples} we apply these
techniques to the following examples: Poisson process with drift
(Section \ref{s:poisson}), Compound Poisson process with Pareto jumps
(Section \ref{s:pareto}), special cases of the ``mixed'' $\alpha$-stable
process (Section \ref{s:astable}), and the generalized inverse Gaussian L\'{e}vy
process (Section \ref{s:gig}).  For each example, we compute one and
two-time moments of $E$.  Methods to compute $U$ for each example
are implemented in MATLAB. The software for doing this is freely
available from the authors and its use is described in
Section \ref{s:guide}.

\section{Background} \label{s:background}

%In the following we consider the Laplace transform, ${\cal{L}}$, which
%is defined for all real valued functions $f$ such that $f(t) = 0$ for
%$t < 0$, and 
%\begin{equation}
%\int_{-\infty}^\infty f(t) e^{-\lambda t} dt < \infty.
%\end{equation}
%This definition can also be extended to Borel measures which have support in
%$[0,\infty)$.  If $\nu$ is such a measure, then the Laplace transform
%of $\nu$ is given by $\int_{-\infty}^\infty e^{-\lambda t} d\nu(t)$
%\footnote{We are avoiding using the limits ``$0$ to $\infty$'' in the
%  definition of the Laplace transform because we will be considering
%  functions which might be discontinuous at $0$.  In these cases,
%  integration by parts is not justified.} .  

We start by recalling the fundamental relationship between a subordinator and
its inverse.  If $D$ is strictly increasing, then for
$t_1,t_2,\dots,t_n$ and $s_1,s_2,\dots,s_n$ positive,
\begin{equation}\label{e:fundeqn}
\{ D(s_i) < t_i, i=1,\dots n \} = \{E(t_i) > s_i , i = 1,\dots,n \}.
\end{equation}
If $D$ is not strictly increasing, then the above relationship holds
off a set of measure $0$ in $(t_1,t_2,\dots,t_n)$ (see 
\cite{veillette:2008}).   

As in the introduction, define the \textit{renewal function} to be the
mean of the inverse subordinator, $U(t) = \mathbb{E}E(t)$ for $t \geq 0$.  Also, let $H_s(t) = P[D(s) < t]$.  From the
L\'{e}vy-Khintchine formula, we have $\int_0^\infty e^{-\lambda t}
dH_s(t) = e^{-s \phi(\lambda)}$.  This together with (\ref{e:fundeqn})
lets us compute the Laplace
transform, $\widetilde{U}$, of $U$:
\begin{eqnarray} 
\widetilde{U}(\lambda) &=& \int_0^\infty U(t) e^{-\lambda t} dt \\
&=& \int_0^\infty \int_0^\infty P[E(t) > s] e^{-\lambda t} ds dt \\
&=& \int_0^\infty \int_0^\infty P[D(s) < t] e^{-\lambda t} ds dt \\
&=& \int_0^\infty \int_0^\infty \frac{1}{\lambda} e^{-\lambda t} dH_s(t)
ds \\
&=& \frac{1}{\lambda} \int_0^\infty e^{-s \phi(\lambda)} ds =
\frac{1}{\lambda \phi(\lambda)}. \label{e:lt_mfpt}
\end{eqnarray}

Thus, $\widetilde{U}$ characterizes the process $E$ (since $\phi$
characterizes $D$).  Since $E$ is non-decreasing, we can define the
Borel measure $dU$ which is induced by $U$, which is commonly referred
to as the \textit{renewal measure}.  Notice that the renewal measure
has the following property.  For a.e. $ 0 \leq a < b$, we have
\begin{eqnarray}
 \int_0^\infty \bone_{(a,b]}(\tau) dU(\tau) &=& U(b) - U(a) \nonumber \\
&=& \int_0^\infty P[E(b) > s] - P[E(a) > s] ds \nonumber \\
&=& \int_0^\infty P[D(s) \leq b] - P[D(s) \leq a] ds \nonumber \\
&=& \mathbb{E} \int_0^\infty \left( \bone_{(-\infty,b]}(D(s)) -
  \bone_{(-\infty,a]}(D(s)) \right)
ds \nonumber \\
&=& \mathbb{E} \int_0^\infty \bone_{(a,b]}(D(s)) 
ds. \label{e:renewalmeasure}
\end{eqnarray} 
By approximating with step functions, this relationship can be
extended to continuous functions $g$ as $\int_0^\infty g(\tau)
dU(\tau) = \mathbb{E}\int_0^\infty g(D(s)) ds$.  With this, we see the Laplace transform of $dU$ is given
by
\begin{equation}
\int_0^\infty e^{-\tau \lambda} dU(\tau) = \mathbb{E} \int_0^\infty e^{-\lambda
  D(s)} ds = \int_0^\infty e^{-s \phi(\lambda)} ds =
\frac{1}{\phi(\lambda)}.
\end{equation}

The pair $U$ and $dU$ can be used to compute all joint moments of the
process $E$.  For non-negative integers $m_1,\dots,m_n$ define 
\begin{equation}\label{e:nmoment}
U(t_1,\dots,t_n;m_1,\dots,m_n) =
\mathbb{E}E(t_1)^{m_1}E(t_2)^{m_2}\dots E(t_n)^{m_n}.
\end{equation}
The order of the moment in (\ref{e:nmoment}) is defined to be $N =
\sum_{i=1}^n m_i$. The following theorem from \cite{veillette:2008}
 gives the $n$-time Laplace transform
 $\widetilde{U}(\lambda_1,\dots,\lambda_n;m_1,\dots,m_n)$ in terms of
 a strictly lower older moment.  

\begin{thm}\label{t:inversemoments}
Let $D$ be a general L\'{e}vy subordinator with L\'{e}vy exponent $\phi$ and let $E$ be the inverse
subordinator of $D$.  The $n$-time Laplace Transform
of the $N^{th}$ order moment  $U(t_1,\dots,t_n;m_1,\dots,m_n)$ defined
in (\ref{e:nmoment}) is given by
\begin{equation}\label{e:inversenmoments}
\widetilde{U}(\lambda_1,\dots,\lambda_n;m_1,\dots,m_n) =
\frac{1}{\phi(\lambda_1+\dots+\lambda_n)} \sum_{i=1}^n m_i
\widetilde{U}(\lambda_1,\dots,\lambda_n;m_1,\dots,m_{i-1},m_i-1,m_{i+1},\dots,m_n).
\end{equation}
\end{thm}

Notice that the Laplace transform of $U(t_1,\dots,t_n;m_1,\dots,m_n)$
is given as the product of $1/\phi$ and the Laplace transform of a
strictly lower order moment.  Taking inverse Laplace transforms, the
$N^{th}$ order moment is given as the sum of convolutions
\begin{equation}
U(t_1,\dots,t_n;m_1,\dots,m_n) = \sum_{i=1}^n m_i \int_0^{\min_i t_i}
U(t_1-\tau,\dots,t_n-\tau;m_1,\dots,m_{i-1},m_i - 1,m_{i+1},\dots,m_n)
dU(\tau).
\end{equation}
For full details, see \cite{veillette:2008}.  

Thus, if one knows the function $U(t;1) = U(t)$, $t \geq 0$, then all
higher order moments can be obtained inductively.  For example, the covariance
is given by
\begin{equation}\label{e:inversecov1}
\mbox{Cov}(E(t_1),E(t_2)) = \int_0^{t_1 \wedge t_2} (U(t_1 - \tau) +
U(t_2 - \tau)) dU(\tau) - U(t_1)U(t_2).
\end{equation}

It is not always easy, however, to invert the Laplace transform
(\ref{e:lt_mfpt}) to obtain $U(t)$
analytically.  In the following, we give two methods for numerically
inverting this Laplace Transform given only the drift $\mu$ and
L\'{e}vy measure $\Pi$ of the subordinator $D$.  Once $U$ is obtained,
the density of the renewal measure can be obtained by numerical
differentiation and integral expressions such as (\ref{e:inversecov1})
can be approximated by numerical integration.  The first method is based on
the Bromwhich integral which expresses the inverse Laplace transform
as a path integral in the complex plane, and the second is based on
the Post-Widder inversion formula, which expresses $U$ as a limit of
terms involving derivatives of $\widetilde{U}$.

\section{Computing $U(t)$}\label{s:numerical}

In Section \ref{s:background} we saw that all moments of an inverse
subordinator can be computed if one has first computed the renewal function
$U(t) = \mathbb{E}E(t)$, which is given by the inverse Laplace transform of
(\ref{e:lt_mfpt}).   In some cases, an
analytical expression for $U$ can be found, and in most cases, the
asymptotics of $U$ can be studied using a Tauberian Theorem.  In
this section we give two methods for calculating $U$ numerically.  The first is
simple and precise, but is difficult to use when $\phi$ is a
complicated function.  The
second is more robust, but requires smoothness in $U$ to be effective. 

\subsection{Method 1: Numerical Integration}\label{s:integration}

The first method involves calculating the inverse Laplace
transform of $\widetilde{U}$ by numerically approximating the
Bromwich integral:
\begin{equation}\label{e:bromwich}
U(t) = \frac{1}{2 \pi i} \int_{b - i \infty}^{b+i\infty} e^{z t}
\widetilde{U}(z) dz,
\end{equation}
where $b$ is chosen such that $\widetilde{U}$ is analytic in the region $\mbox{Re}(z) \geq b$.  Using
the fact that $U(t) = 0$ for $t < 0$, and symmetry properties of analytic functions which are real valued
on the real axis,  (\ref{e:bromwich}) simplifies into two equivalent expressions:
\begin{equation}\label{e:integral}
U(t) = \frac{2 e^{bt}}{\pi} \int_0^\infty \left( \mbox{Re}(
  \widetilde{U}(b+iu)) \cos(u t) \right) du = \frac{-2 e^{bt}}{\pi}
\int_0^\infty  \left( \mbox{Im}(
 \widetilde{U}(b+iu)) \sin(u t) \right) du.
\end{equation}
For details see \cite{Grassmann:2000}.

 From the
L\'{e}vy Khintchine formula (\ref{e:LKformula}), we have
\begin{eqnarray}
\phi(b + iu) &=& \mu b + i \mu u + \int_0^\infty \left( 1 - e^{-x(b + i
    u)} \right) \Pi(dx) \\
&=& \left( \mu b + \int_0^\infty\left( 1 - e^{-x b} \cos( x u) \right) \Pi(dx)
\right) + i \left(\mu u - \int_0^\infty e^{-x b} \sin( x u) \Pi(dx)
\right)  \label{e:phirphii} \\
& \equiv & \phi_r(b+i u) + i \phi_i(b + i u). 
\end{eqnarray}  
Now, since $\widetilde{U}(\lambda) = (\lambda \phi(\lambda))^{-1}$,
a simple calculation gives
\begin{eqnarray}
\mbox{Re}( \widetilde{U}(b + i u)) &=& \frac{b \phi_r(b + i u) - u
  \phi_i(b + i u)}{(b^2 + u^2)(\phi_r(b+i u)^2 + \phi_i(b + i u)^2)} \label{e:realut} \\
\mbox{Im} (\widetilde{U}(b + i u) )&=& \frac{b \phi_i(b + i u) + u
  \phi_r(b + i u)}{(b^2 + u^2)(\phi_r(b+i u)^2 + \phi_i(b + i u)^2)}. \label{e:imut}
\end{eqnarray}

Since $\phi(0) = 0$ and $\phi$ is increasing in $\lambda$,
$\widetilde{U}(\lambda) = (\lambda \phi(\lambda))^{-1}$ has a
singularity the origin.  Therefore, we much choose $b > 0$ above.
Then, given $\mu$ and $\Pi$, we evaluate the integrals in (\ref{e:phirphii}) to
obtain $\phi_r$ and $\phi_i$ and use (\ref{e:realut}) and (\ref{e:imut}) to obtain
$\mbox{Re}(\widetilde{U}(b + iu))$ and $\mbox{Im}(\widetilde{U}(b+iu))$
for fixed $u$ and then compute either integral in (\ref{e:integral}) to get
$U(t)$. Alternatively, if $\phi$ is known in closed form, then $\phi_r$
and $\phi_i$ can be computed directly.  This method works fairly well when $\phi_r$ and $\phi_i$ are
easy to compute, for instance, in the case of the Poisson
process with drift.  The main problem with the method is when the integrands in (\ref{e:integral})
and/or (\ref{e:phirphii}) are
highly oscillatory and slowly decaying, causing most integration
algorithms to converge very slowly.

\subsection{Method 2: Post-Widder Inversion}\label{s:postwidder}
The following method is based on the Post-Widder inversion formula
(\cite{Grassmann:2000}, Theorem 2 or ~\cite{feller:1971}, section
VII.6). In order to justify using this formula in our case, we state and
prove this result under slightly weaker conditions than those found in
~\cite{feller:1971}.     
 
\begin{thm}
(Post-Widder Inversion).  Let $u: \mathbb{R}^+
\rightarrow \mathbb{R}^+$ be a continuous function such that $u(x)/x
\leq C$ as $x \rightarrow \infty$ for some $C \geq 0$, and let $\widetilde{u}(\lambda)$ be the Laplace
transform of $u$.  Then for $t>0$, we have
\begin{equation}\label{e:pwinversion}
u(t) = \lim_{k \rightarrow \infty} \frac{(-1)^{k-1}}{(k-1)!} \left(
  \frac{k}{t} \right) ^{k} \widetilde{u}^{(k-1)}\left( \frac{k}{t}
\right),
\end{equation}   
where  $\widetilde{u}^{(k)}$ denotes the $k^{th}$ derivative, $k \in \mathbb{Z}^+$.  
\end{thm}

\begin{pf}
Fix $t>0$.  Let $X_i$, $i=1,2,\dots$ be iid gamma random variables
with mean $t$ and variance $t^2$ and let $M_k =
\frac{1}{k} \sum_{i=1}^k X_i$.  The idea is to approximate $u(t)$ by
$\mathbb{E}u(M_k)$.  Notice that $M_k$ has a gamma
distribution with mean $t$ and variance $t^2/k$.  By the
law of large numbers,  $M_k \rightarrow^d t$ as $k \rightarrow \infty$ and by the continuous
mapping theorem (\cite{gut:2005}, Theorem 5.10.4),  $u(M_k)
\rightarrow^d u(t)$.  Observe that we also have the following
convergence in mean: 
\begin{equation}\label{e:cmt}
\mathbb{E}u(M_k) \rightarrow u(t),
\end{equation}
which follows from the uniform integrability of $u(M_k)$.  This can be
seen using the
density of a gamma random variable and the assumption that $u(x) \leq Cx$
as $x \rightarrow \infty$.  Indeed, we have
\begin{eqnarray}
\lim_{a \rightarrow \infty} \sup_{k \geq 1} \mathbb{E} u(M_k)
\bone_{M_k > a} &=& \lim_{a \rightarrow \infty} \sup_{k \geq 1}
\int_a^\infty u(x) x^{k-1} e^{-k x/t} dx \\
& \leq &  \lim_{a \rightarrow \infty} \sup_{k \geq 1} C \int_a^\infty
x^k e^{-k x/t} dx \\
&=& \lim_{a \rightarrow \infty} \sup_{k \geq 1} C \int_a^\infty
(x e^{-x/t})^k dx \\
&=&  \lim_{a \rightarrow \infty}  C \int_a^\infty
(x e^{-x/t}) dx  \\
&=&  0.
\end{eqnarray}

Now, (\ref{e:cmt}) implies 
\begin{eqnarray}
u(t) = \lim_{k\rightarrow \infty}\mathbb{E}u(M_k) &=& \lim_{k\rightarrow \infty} \frac{1}{(k-1)!} \left(\frac{k}{t} \right)^k
\int_0^\infty u(x) x^{k-1} e^{-k x/t} dx \\
&=& \lim_{k\rightarrow \infty} \frac{1}{(k-1)!} \left(\frac{k}{t} \right)^k \left( (-1)^{k-1}
  \frac{\partial^{k-1}}{\partial \lambda^{k-1}} \int_0^\infty u(x)
  e^{-\lambda x} dx \right) \Big |_{\lambda = k/t} \\
&=& \lim_{k\rightarrow \infty} \frac{(-1)^{k-1}}{(k-1)!} \left(\frac{k}{t} \right)^k
\widetilde{u}^{(k-1)}(k/t).
\end{eqnarray}
This verifies (\ref{e:pwinversion}).  
\end{pf}

\begin{rems} 
 \textbf{1.}  In order to justify using the Post-Widder formula, we must be
sure that the renewal function $U$ satisfies $U(t)/t \leq C$ as $t
\rightarrow \infty$.
Fortunately, the Renewal Theorem (see for instance, Proposition 3 in
\cite{Nord:2005}) implies that if $D$ has finite mean, then $U(t)$
will grow as $O(t)$ as $t \rightarrow \infty$.  If $D$ has
  infinite mean, its mean first hitting time $\mathbb{E}E(t)$ will, on
average, grow slower than in the finite mean case for large $t$.  To
see this, any subordinator $D$ with infinite mean can be bounded below
by a subordinator $\hat{D}$ with finite mean by truncating the large
jumps of $D$\footnote{That is, $\hat{D}(s) \leq D(s)$ for all $s \geq
  0$ a.s.}.  The inverse, $\hat{E}$, of $\hat{D}$ will be a.s. greater
than $E$, and thus $U(t) \leq \hat{U}(t) = O(t)$.

\textbf{2.}  One must check that $U$ is
continuous before using the Post-Widder formula.  Proposition A.1 in
\cite{veillette:2008} implies that $U$ is continuous if $D$ is
strictly increasing.     
\end{rems}

Thus,
calculating $U(t)$ involves two steps:

\medskip

\noindent\textbf{1.}  Calculating 
\begin{equation}\label{e:Uks}
U_{k_i}(t) =  \frac{(-1)^{k-1}}{(k_i-1)!} \left(
  \frac{k_i}{t} \right) ^{k_i} \widetilde{U}^{(k_i-1)}\left( \frac{k_i}{t}
\right)
\end{equation}
 for some set of integers $k_1,k_2,\dots,k_n$ (the choice of the $k_i$'s will
 be addressed later).

\medskip
 
\noindent \textbf{2.} Using the values of $U_{k_i}(t)$, $i=1,\dots n$, to approximate the limit in (\ref{e:pwinversion}).

\medskip

  We will first focus on step $1$ above.  The main difficulty with using this method resides in the fact that (\ref{e:pwinversion}) involves
derivatives of arbitrarily high orders, which are often difficult to
compute.  Fortunately, in our case, this calculation can be done in a
reasonable fashion, which we now outline.  Recall
Leibnitz's formula, which states that for smooth functions $f,g$,
\begin{equation}
\frac{d^k}{dt^k} f(t) g(t) = \sum_{i=0}^k {k \choose i}
g^{(k-i)}(t)f^{(i)}(t).
\end{equation}
Let $k \geq 1$,$t>0$, and $\psi(\lambda) = 1/\phi(\lambda)$.
Applying Leibnitz's formula to $\widetilde{U}(\lambda) =
\lambda^{-1} \psi(\lambda)$ with $\lambda  = k/t$ gives
\begin{eqnarray}
 \frac{(-1)^{k-1}}{(k-1)!} \left(
  \frac{k}{t} \right) ^{k} \widetilde{U}^{(k-1)}(k/t) &=&  \frac{(-1)^{k-1}}{(k-1)!} \left(
  \frac{k}{t} \right) ^{k} \sum_{i=0}^{k-1}  {k-1 \choose i} \left(
\frac{ (-1)^{k-1-i} (k-1-i)! }{(k/t)^{k-i}}
\right) \psi^{(i)}(k/t)  \\
&=&  \sum_{i=0}^{k-1} (-1)^{i}  \frac{k^i }{i! t^i} \psi^{(i)}(k/t) \\
&=& \mbox{\textbf{v}}_k \cdot \mbox{\textbf{w}}_k, \label{e:dotproduct}
\end{eqnarray}
where the vectors $\textbf{v}_k,\textbf{w}_k \in \mathbb{R}^k$ are
given by
 \begin{eqnarray}
({\textbf{v}}_k)_i &=&   \left(    \frac{  (-1)^{i} k^i}{i! t^i } 
\right), \quad i = 0,\dots, k-1   \label{e:vk} \\
({\textbf{w}}_k)_i &=& \left( \psi^{(i)}(k/t) \right), \quad i =0,\dots,k-1.
\end{eqnarray}

To compute the components of $\textbf{w}_k$, we use an idea from
  \cite{Leslie:1991}.  Using Leibnitz's formula on the unit function
  $1 = \phi(\lambda)\psi(\lambda)$, we
  see that for any $\lambda >0$,
\begin{eqnarray} 
\sum_{i=0}^j { j \choose i} \phi^{(j-i)}(\lambda)
\psi^{(i)}(\lambda) &=& \frac{\partial^{j}}{\partial \lambda^{j}} \phi(\lambda) 
\psi(\lambda)  \\
&=&  \frac{\partial^{j}}{\partial \lambda^{j}} 1 \\
&=& \begin{cases} 1 \quad j=0, \\
0 \quad j \geq 1 \end{cases}, \quad j = 0,\dots,k-1.
\end{eqnarray}  
From this, we obtain the following matrix equation:
\begin{equation}\label{e:matrixeqn}
\left(
\begin{array}{cccccc}
\phi(\lambda) & 0 & 0 & 0 & \dots & 0 \\
\phi'(\lambda) & \phi(\lambda) & 0 & 0 & \dots & 0 \\
\phi''(\lambda) & 2 \phi'(\lambda) & \phi(\lambda) & 0 & \dots & 0 \\
\phi'''(\lambda) & 3 \phi''(\lambda) & 3 \phi'(\lambda) &
\phi(\lambda) & \dots & 0 \\
\vdots & & & & \ddots & \\
\phi^{(k-1)}(\lambda) & (k-1) \phi^{(k-2)}(\lambda) &  &\dots& & \phi(\lambda) \end{array}
\right) \left( \begin{array}{c} \psi ( \lambda ) \\
 \psi' ( \lambda ) \\
 \psi'' ( \lambda ) \\
 \psi''' ( \lambda )  \\
\vdots \\
 \psi^{(k-1)}(\lambda)  \end{array} \right) =  \left( \begin{array}{c} 1 \\
 0 \\
 0 \\
 0  \\
\vdots \\
 0 \end{array} \right).
\end{equation}
Choosing $\lambda = k/t$, we conclude that $\textbf{w}_k$ satisfies the
matrix equation $G_k \textbf{w}_k = e_k$, where $e_k = (1,0,\dots,0)'
\in \mathbb{R}^k$ and
\begin{equation}\label{e:gji}
(G_k)_{ji} = \begin{cases} {j-1 \choose i-1} \phi^{(j-i)}(k/t), \quad
  1 \leq i \leq j \leq k \\
0, \quad \quad 1 \leq j < i \leq k \end{cases}.
\end{equation}
Observe that the entries in $G_k$ are easily expressed in terms of
the drift $\mu$ and the L\'{e}vy measure $\Pi$ corresponding to the subordinator
$D$.  Indeed, we have
\begin{equation} \label{e:phidiff}
\phi^{(j)}(k/t) = \begin{cases} \displaystyle \mu \frac{k}{t} + \int_0^\infty \left(1 - e^{-kx/t}
\right) \Pi(dx), \quad j=0, \\
\\
\displaystyle \mu + \int_0^\infty x e^{-kx/t} \Pi(dx), \quad j=1, \\
\\
 \displaystyle (-1)^{j+1} \int_0^\infty x^j e^{-kx/t} \Pi(dx),\quad j \geq 2 \end{cases}.
\end{equation}
Thus, to compute $\textbf{w}_k$, one needs to compute the
entries of $G_k$ using (\ref{e:gji}) and (\ref{e:phidiff}),
and then solve the matrix equation $G_k \textbf{w}_k = e_k$.  Notice
that the integrals in (\ref{e:phidiff}) are much easier to compute
than those in Section \ref{s:integration} since the integrands do not
oscillate; they decay
exponentially and are positive.  Observe that $\textbf{v}_k$, and hence
$\textbf{v}_k \cdot \textbf{w}_k$ are easy to compute.  Using this
method, we get $U_{k_i}(t) = \textbf{v}_{k_i} \cdot \textbf{w}_{k_i}$
in (\ref{e:Uks}).

\begin{rem}
The density $U'(t)$ of the renewal measure is also of interest, for
example in (\ref{e:inversecov1}).  It can be approximated as $U$ is in
(\ref{e:Uks}),  since the Laplace transform of $U'(t)$ is
$\psi(\lambda)$.    One needs to compute $\psi^{(k-1)}(k/t)$,
which is obtained with no extra cost from (\ref{e:matrixeqn}) since it
is the last component of $\textbf{w}_k$.     
\end{rem}

 For step 2 above, we refer to the
  technique explained in \cite{frolov:1998}.  To summarize,
 it has been shown (\cite{abate:1992},\cite{jagerman:1982}) that if $U$ is smooth, then we have the
 following series expansion:
\begin{equation}\label{e:seriesexp}
U_k(t) = U(t) + \sum_{m=1}^\infty \frac{a_m(t)}{k^m},
\end{equation}
where $a_m(t)$ are remainder terms.  Write $h_i = 1/k_i$ and let
$\widehat{U}_{h_i}(t) = U_{k_i}(t)$.  With this, (\ref{e:seriesexp}) implies
\begin{equation}
\widehat{U}_{h_i}(t) = U(t) + \sum_{m=1}^\infty a_m(t)h_i^m.
\end{equation}
Our goal is to compute $\widehat{U}_0(t) = U(t) = \lim_{k \rightarrow
  \infty} U_k(t)$ with $t$ fixed.  To
do so, we consider $\widehat{U}_{h}(t)$ as a function of $h$ and,
given $h_1,h_2,\dots,h_n$ and
$\widehat{U}_{h_1},\dots,\widehat{U}_{h_n}$, we write down the so-called ``Lagrange polynomial''
$P_n(h)$:  this is the polynomial of degree $n-1$ which passes through the $n$ points
$(h_i,\widehat{U}_{h_i}(t))$, $i=1,2,\dots,n$,
\begin{equation}\label{e:lagrangepoly}
P_n(h) = \sum_{i=1}^n \left( \prod_{j\neq i } \frac{h-h_j}{h_i-h_j}
\right) \widehat{U}_{h_i}(t).
\end{equation}
Observe that if $h = h_k$ for some $k=1,2,\dots,n$, then all the
summands in (\ref{e:lagrangepoly}) vanish except for the term where $i=k$,
which is equal to $\widehat{U}_{h_k}(t)$. Thus the polynomial $P_n(h)$ passes
through the $n$ points $(h_i,\widehat{U}_{h_i}(t))$, $i=1,2,\dots,n$. 

Observe that we cannot compute $\hat{U}_0(t)$ since $h=0$ corresponds
to $k = \infty$.  However, with the method described above, we can compute $\widehat{U}_{h_i}$
for $h_i \neq 0$ and then approximate $\widehat{U}_0(t)$ by $P_n(0)$.
Since $\widehat{U}_0(t) = U(t)$, $U$ is then approximated by the linear combination    
\begin{equation}\label{e:lincomb}
U(t) \approx P_n(0) = \sum_{i=1}^n c^{(n)}_i U_{k_i}(t),
\end{equation}
where the weights $c^{(n)}_i$, calculated by setting $h=0$ in
(\ref{e:lagrangepoly}), depend on the choice of
$\{h_i\}_{i=1}^n$\footnote{As mentioned in \cite{frolov:1998}, the linear
  combination (\ref{e:lincomb}) obtained by polynomial interpolation
  is equivalent to taking the linear combination which cancels the
  first $n$ remainder terms in (\ref{e:seriesexp}).}.

Because
the above method allows us to calculate $\widetilde{U}^{(k)}$ for
large $k$, we shall take
$k_i = 2^{i-1}$, $i=1,2,\dots,n$. The corresponding weights $c_i^{(n)}$
are given by 
\begin{equation} 
c^{(n)}_i = \frac{ (-1)^{n-i} 2^{i(i-1)/2}}{\prod_{j=1}^{i-1} (2^j - 1)
  \prod_{j=1}^{n-i} (2^j-1)}, \quad i = 1,2,\dots,n.
\end{equation}

Thus, given a desired level of accuracy $\epsilon > 0$ and $t>0$ fixed, we
compute the sequence $P_1(0),P_2(0),\dots$ until $|P_n(0) - P_{n-1}(0)| <
\epsilon$ and then set $U(t) = P_n(0)$.  Typically this method
converges with $n < 9$ for reasonable values of $\epsilon$.  Since $k_{10} = 512$, it is often impossible to compute
$U_{k_i}$ for  $i \geq 10$ in double precision arithmetic, as noted in
the remarks below. 

Note that for (\ref{e:seriesexp}) to hold and for this method to be most effective, $U$
must be sufficiently smooth.  We found that for points at which $U$ is
not differentiable, the sequence $P_i$ had a slower rate of convergence.

We close this section with two remarks regarding the computation of $U_k(t)$.    

\begin{rem}
Caution should be taken when computing terms like $k^j/j!$ for $k,j$
large, since numbers like $100^{100}$ and $200!$ are outside the range
of double precision arithmetic.  Instead, one should use expressions
like $k^j/j! =
\exp(j \log(k) - \sum_{i=1}^j \log(i) )$, because while $k^j$ and $j!$
may be large, their ratio may be small. Nevertheless, for large
enough $k$ and $j$ (relation  (\ref{e:vk}), for example, which requires
computing $k^{k-1}/(k-1)!$), even the ratio may be too large.  This
happens for instance if $k = k_{11} = 2^{10} = 1024$, in which case
$1024^{1023}/1023! \approx 6.5 \times 10^{442}$.  To be safe, we used
at most $k_{10} = 512$.
\end{rem}

\medskip

\begin{rem} A similar overflow problem can occur for $k/t$ large.  To avoid this, one can make the following
adjustment to (\ref{e:pwinversion}).  Let $c>0$, then
\begin{eqnarray}
\frac{(-1)^{k-1}}{(k-1)!} \left(
  \frac{k}{t} \right) ^{k} \widetilde{U}^{(k-1)}\left( \frac{k}{t}
\right) &=& \frac{(-1)^{k-1}}{(k-1)!} \left(
  \frac{k}{c t} \right) ^{k} c^{k} \widetilde{U}^{(k-1)}\left( \frac{k}{t}
\right) \\
&=& \frac{(-1)^{k-1}}{(k-1)!} \left(
  \frac{k}{ c t} \right) ^{k} \widetilde{U}^{(k-1)}\left(
  \frac{k}{c t}; c
\right) \label{e:pwinversion2},
\end{eqnarray}
where $\widetilde{U}(\lambda;c)$ is a rescaled version of $\widetilde{U}$:
\begin{equation}
\widetilde{U}(\lambda;c) \equiv c \widetilde{U}(c
\lambda) = \frac{1}{\lambda \phi(c \lambda)}.
\end{equation}
Repeating the steps of this method with (\ref{e:pwinversion2}), we get
\begin{equation}
U_k(t) =  \widetilde{\textbf{v}}_k \cdot
\widetilde{\textbf{w}}_k,
\end{equation}
where $(\widetilde{\textbf{v}}_k)_i = (\textbf{v}_k)_i/c^i$ and
$\widetilde{\textbf{w}}_k$ is the solution to the matrix equation
$\widetilde{G}_k \widetilde{\textbf{w}}_k = e_k$, where
$(\widetilde{G}_k)_{ji} = (G_k)_{ji}c^{j-i}$.  To compute these
products and ratios, a procedure as the one describes in the previous
remark should be done to avoid overflow.  We found
the choice $ c = t^{-1}$ useful.  

\end{rem}

\subsection{Testing the methods}
  
 In Table
\ref{t:errors}, the Post-Widder method is tested in cases where $U$ can be
computed explicitly (see Section \ref{s:examples} for explanations of
each example).     The values in the table were obtained using a threshold of $\epsilon
  = 10^{-8}$.  With the exception of starred entry in Table \ref{t:errors}, all converged
  within this threshold with $n \leq 9$. If sufficient convergence did not occur,
  $P_9(t)$ is used as the approximation.  In these cases, the approximation is very good.

\begin{figure}[ht]
\centering
\includegraphics[width=0.8\textwidth]{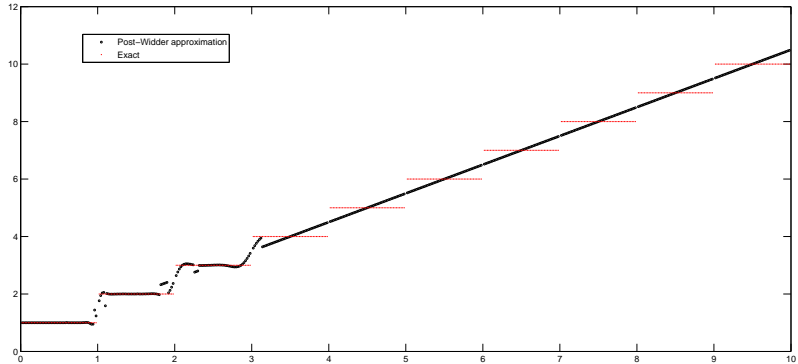}
\caption{A comparison of the true value of $U$ for a Poisson process
  with no drift and its approximation given by the Post-Widder
  method.  Since $U$ in this case has little regularity, the quality
  of the approximation is not always good. }
\label{f:pwcompare}
\end{figure}

  For the  Poisson process, $U$ is discontinuous at the
  integers. Therefore, we
  used instead $t=1.1,10.1,100.1$.  Since the Post-Widder method does not apply
  when $U$ is discontinuous, we do not expect to get good results in
  the Poisson case.  And indeed, the absolute errors are large.  For
  instance, Table \ref{t:errorsint} shows that when $t=10.1$, the algorithm converged but to the wrong
  value (the absolute error between the true and computed value is
  $0.4$).   The plot in Figure \ref{f:pwcompare} shows in more detail the erratic
  behavior of the Post-Widder method in this case.  Notice that for
  larger times this method approximates $U$ (which in this case is a
  step function) with a straight line.

\begin{table}[ht]
\centering   
  \begin{tabular}{| l | c | c | c | c | c | }
    \hline
    Subordinator & $t = 0.01$ & $t = 0.1$ & $t = 1$ & $t = 10$ & $t =
    100$  \\ \hline 
    $\alpha$-stable ($\alpha = 0.5$) & $2.5 \times 10^{-13}$ & $8.2 \times 10^{-13}$ &
    $1.7 \times 10^{-13}$ & $7.4 \times 10^{-13}$ & $2.6 \times
    10^{-11}$ \\ \hline
    Uniform Mixture of $\alpha$-stable & $7.8 \times 10^{-12}$ & $6.6 \times 10^{-11}$ & $1.4
    \times 10^{-10}$ & $1.5 \times 10^{-9}$ & $8.7 \times 10^{-10}$ \\
    \hline
    Gamma Process ($\kappa = \gamma = 1$) & $1.7 \times 10^{-13}$  &
    $7.0 \times 10^{-14}$ & $1.2 \times 10^{-12}$ & $3.4 \times
    10^{-10}$ & $9.6 \times 10^{-12}$ \\ \hline
    Inverse Gaussian($\gamma=\delta = 1$) &$1.3 \times 10^{-5*}$ &
    $1.7 \times 10^{-12}$ & $1.3 \times 10^{-11}$ & $6.7 \times
    10^{-11}$ & $3.0 \times 10^{-9}$ \\ \hline
 \hline
  \end{tabular}
\caption{Absolute errors between exact values of $U$ and those given
  by this method for a selection of subordinators whose mean first
  passage time can be computed exactly.  Here we used a threshold of
  $\epsilon = 10^{-8}$.   }
\label{t:errors}
\end{table}

In the cases where $U$ has little regularity,
numerical integration (method 1) is the more accurate method.  To see
this, we applied the numerical integration method (using the MATLAB
function \texttt{quad}) to the case of the
Poisson process and obtained the absolute errors in Table
\ref{t:errorsint}.  To compute these, we approximated the oscillatory integral in
(\ref{e:integral}) by partitioning the range of integration $(0,\infty)$
into intervals of the form $I_0 = (0, \pi/(2t)$, $I_k = [k \pi/(2 t),
(k+2)\pi/(2t)]$, $k \geq 1$, $k$ odd, and then
summing the integrals over each $I_k$ until the contribution 
over one such interval
became less then a threshold $\epsilon > 0$.  Due to the slowly decaying
nature of the integrand, obtaining convergence for $\epsilon <
10^{-6}$ is difficult.

  We've found that integration is not feasible in many other
cases (except possibly the $\alpha$-stable case). 

\begin{table}[ht]
\centering   
  \begin{tabular}{| l | c | c | c | c | c | }
    \hline
    Subordinator & $t = 0.01$ & $t = 0.1$ & $t = 1.1$ & $t = 10.1$ & $t =
    100.1$  \\ \hline 
    Poisson Process (integration)  & $  1.08 \times 10^{-6}  $ & $ 1.57
    \times 10^{-6}  $ &
    $ 3.6 \times 10^{-11} $ & $ 4.25 \times 10^{-6} $ & $ 1.09 \times
    10^{-4} $ \\ \hline
    Poisson Process (Post-Widder) & $ < 10^{-16}$ & $7.1 \times 10^{-10}$ &
    $3.7 \times 10^{-2 *}$  & $0.4$& $0.4$ \\ \hline

  \end{tabular}
\caption{Absolute errors between exact values of $U$ and those given
  by numerical integration for the poisson process (whose mean first
  passage time can be computed exactly).  Here we used a threshold of
  $\epsilon = 10^{-6}$.   }
\label{t:errorsint}
\end{table}

\section{Obtaining $U(t) = \mathbb{E}E(t)$ and $\mathrm{Corr}(E(s),E(t))$ for various inverse
  L\'{e}vy subordinators}\label{s:examples}

We now present various examples of L\'{e}vy subordinators $\{ D(s), s
\geq 0\}$ and
their inverses $\{ E(t), t \geq 0 \}$ and  calculate the one time
moment $U(t) = \mathbb{E}E(t)$ and the correlation function
$\mathrm{corr}(E(s),E(t))$.  We found the table of Laplace transforms
\cite{Roberts:1966} useful for some of the following calculations.

Since we will encounter Laplace transforms which cannot be inverted
analytically, we will study their asymptotics using a Tauberian
Theorem (\cite{Bertoin:1996} page 10), which we state here for convenience.  Recall that a function
$\ell(t)$, $t>0$, is \textit{slowly varying} at $0$ (respectively
$\infty$) if for all $c>0$, $\lim(\ell(ct)/\ell(t)) = 1$ as
$t\rightarrow 0$ (respectively $t\rightarrow \infty$). 
\begin{thm} 
 (Tauberian Theorem) Let $\ell:(0,\infty) \rightarrow (0,\infty)$ be a slowly varying
function at $0$ (respectively $\infty$) and let $\rho \geq 0$.  Then
for a function $U: (0,\infty) \rightarrow (0,\infty)$, the
following are equivalent:

(i) $U(x) \sim x^\rho \ell(x)/\Gamma(1 + \rho), \quad x \rightarrow 0$
(respectively $x \rightarrow \infty$).

(ii) $\widetilde{U}(\lambda) \sim \lambda^{-\rho-1} \ell(1/\lambda), \quad \lambda \rightarrow \infty$
(respectively $\lambda \rightarrow 0$).
\end{thm}

\subsection{Poisson process}\label{s:poisson}

Consider the process $D(s) = \mu s + N(s)$, where $\{N(s),
s \geq 0 \}$ is a Poisson process with rate $r>0 $.  The L\'{e}vy exponent for this is given by 
\begin{equation}\label{e:phiforpois}
\phi(\lambda) = \mu \lambda +  r(1 - e^{-\lambda}),
\end{equation}
and the L\'{e}vy measure $\Pi(dx)$ is a point-mass at $x=1$ with
weight $r$.  The sample paths of $D$ are straight lines with
slope $\mu$, together with jumps of size $1$ which happen at random
times $\{ \tau_k\}_{k=1}^\infty $, such that 
$\{\tau_i-\tau_{i-1}\}_{i=1}^\infty$ are iid exponential with mean
$1/r$.   

\begin{figure}[ht]
\centering
\includegraphics[width=0.8\textwidth]{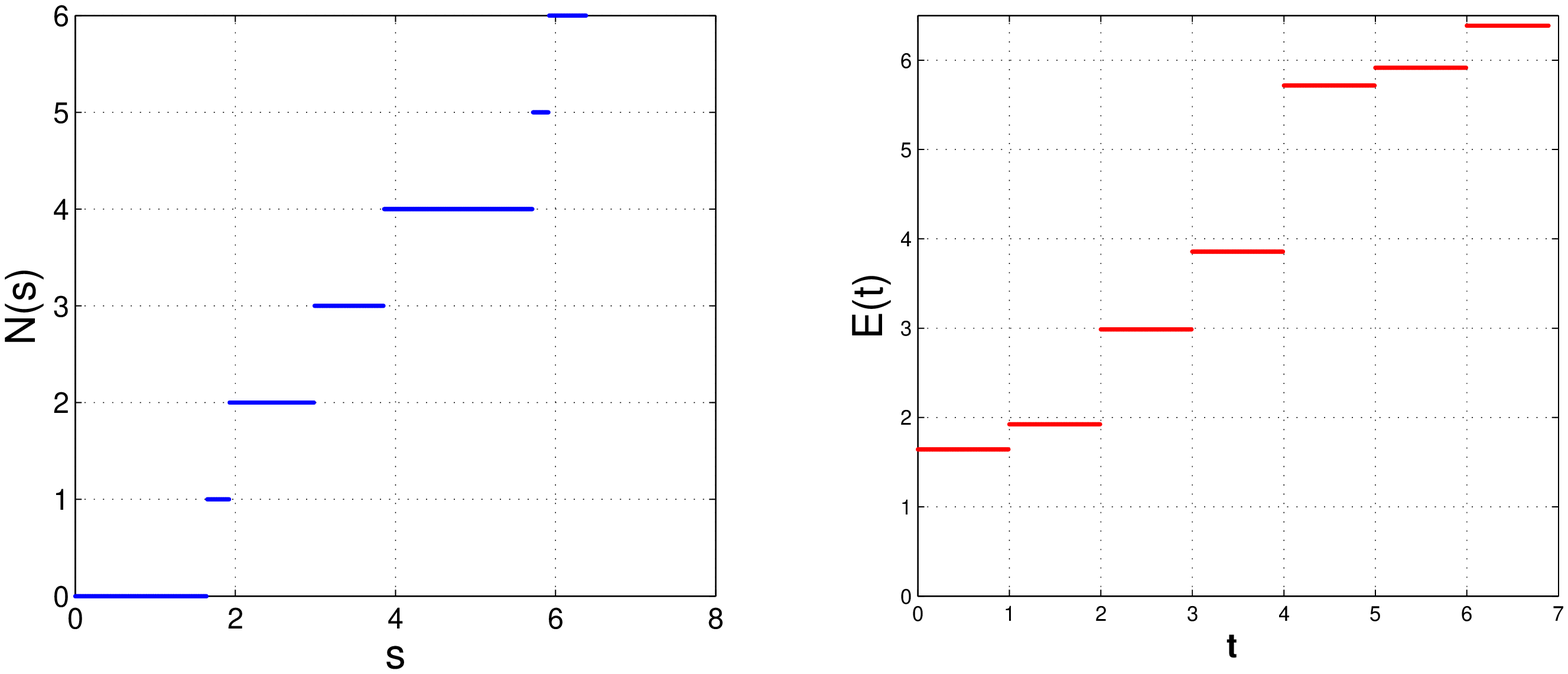}
\caption{A sample path of the Poisson process $N(s)$ together with its
inverse.}
\label{f:poisplot}
\end{figure}

First consider the case with no drift, $\mu =0$. Figure
\ref{f:poisplot} displays a sample path of the Poisson process $N(s)$
together with its inverse $E(t)$ which is the first time $N(s)$
exceeds the level $t$.  Notice that each segment in the
plot of $E$ has length $1$.  For any $0<t_0<1$, one
must wait a random time $\tau_1$ for the process $N$ to surpass $t_0$,
thus $E(t_0) = \tau_1$.  More generally, the
inverse subordinator $E(t)$ is a sum of $\lfloor t+1 \rfloor$ iid
exponential random variables with mean $1/r$, implying $ E(t) \sim
\Gamma(\lfloor t+1 \rfloor,1/r)$.

Using the density of a Gamma random
variable, the $\gamma$-moment of the $E(t)$, with $\gamma >0$, is given by
\begin{eqnarray}
\mathbb{E}E(t)^\gamma &=& \frac{r^{\lfloor t+1 \rfloor}}{\Gamma(\lfloor t+1 \rfloor) } \int_0^\infty x^\gamma x^{ \lfloor t+1 \rfloor - 1} e^{-x r} dx  \\
&=& \frac{ r^{\lfloor t+1 \rfloor} \Gamma(\lfloor t+1 \rfloor+
  \gamma)}{\Gamma(\lfloor t+1 \rfloor) r^{\lfloor t+1 \rfloor + \gamma} } \\
&=& \frac{\Gamma(\lfloor t+1 \rfloor+ \gamma) }{ r^ \gamma \Gamma(\lfloor t+1 \rfloor)} .
\end{eqnarray}
Setting $\gamma = 1$ yields $U(t) = \lfloor t+1 \rfloor/r$ and the renewal measure, $dU(t)$, for the Poisson process is the
measure which assigns a mass of $1/r$ to each integer $n \geq 0$.
The covariance, (\ref{e:inversecov1}), of this process is then given by
\begin{equation}\label{e:covnodrift}
\mbox{\rm Cov}(E(s),E(t)) = \left( \sum_{k=0}^{\lfloor s \wedge t \rfloor} \frac{U(s-k)}{r} +
  \frac{U(t-k)}{r} \right)  - U(s) U(t).
\end{equation}

\begin{figure}[ht]
\centering
\includegraphics[width=1.1\textwidth]{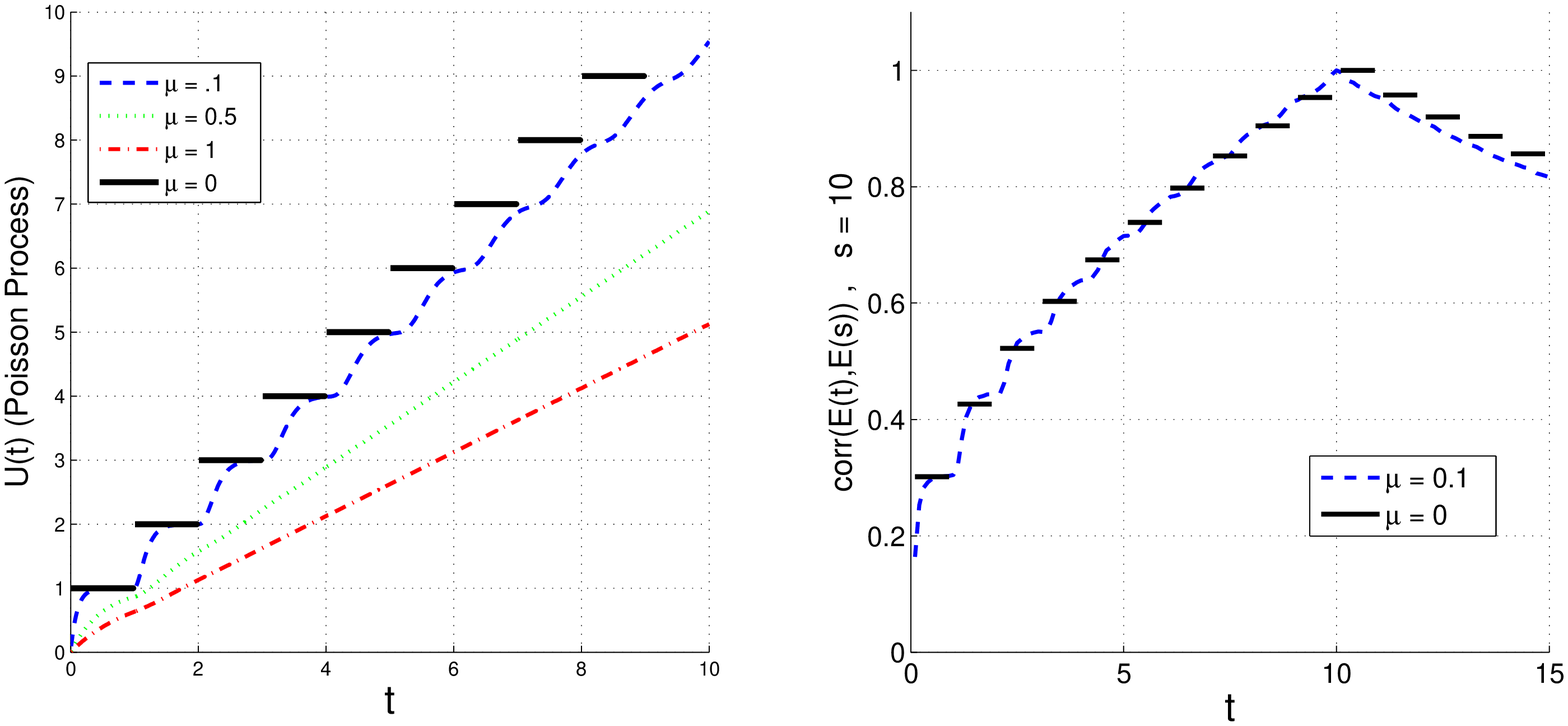}
\caption{(Left) Plots of $\mathbb{E}E(t)$, where $E$ is the inverse of a
  Poisson process with rate $r=1$ and various drifts.  For $\mu >0$,
  these plot were generated using the numerical integration method in
  Section \ref{s:integration}.  (Right)  The correlation coefficient $\mbox{corr}(E(s),E(t))$
  of the Poisson process with $s=10$ fixed and $t$ varying with drifts $\mu =0$ and $\mu = 0.1$.  This
  was calculated using (\ref{e:covnodrift}) in the case $\mu =0$ and
  by numerical approximation of the integral in equation
  (\ref{e:inversecov1}) for $\mu >0$.  Since here $s = 10$, the
  correlation equals $1$ when $t=10$ as well.}
\label{f:poismfpt}
\end{figure}

Now assume positive drift, i.e. $\mu >0$.  From
(\ref{e:lt_mfpt}) and (\ref{e:phiforpois}), the Laplace transform of $U(t)$ is given by
\begin{equation}\label{e:poisfpt}
\widetilde{U}(\lambda) = \frac{1}{\lambda( \mu \lambda + r(1- e^{-\lambda}))}.
\end{equation}
Due to the simple form of $\widetilde{U}(t)$, $U(t)$
can be calculated by numerical integration (i.e. method 1 in Section
\ref{s:numerical}) \footnote{The Post-Widder method also works well
  when $\mu/r \gg 0$.  It fails to converge however, when $\mu/r$ is
  close to $0$.}.  In Figure
\ref{f:poismfpt} we plot $U(t) = \mathbb{E}E(t)$ for various values of
$\mu$, namely $\mu = 0,0.1,0.5,1$.  We also compute the correlation
coefficient.  It is obtained using (\ref{e:covnodrift}) in the case of
no drift but it is not known in closed for $\mu > 0$.  To compute for
$\mu >0$, we use (\ref{e:inversecov1}).  The integral in (\ref{e:inversecov1}) involves the
renewal measure $dU$, but since in this case the function $U(t)$ is strictly increasing and continuous for $\mu$
positive, $dU$ is given by
$U'(t)dt$.  We obtain $U'(t)$ here by discrete approximation.  The correlation coefficient
$\mbox{corr}(E(s),E(t))$ is plotted in Figure \ref{f:poismfpt} as a
function of $t$ with $s=10$ fixed for
zero and non-zero drift.  

%\begin{figure}[ht]
%\centering
%\includegraphics[width=0.8\textwidth]{poisoon_cov.eps}
%\caption{ (a) $U'(t)$ for
%  the Poisson process with drifts $\mu = .1,.5$. (b) Plot of $r(s,t)$
%  with $s = 10$ and $t$ varying. From top to bottom corresponds to $\mu = 0,.1,%.5$.  }
%\label{f:poiscov}
%\end{figure}

\subsection{Compound Poisson processes}  \label{s:pareto}

Let $\xi_i$ $i=1\dots$ be positive iid random variables with probability measure
$\nu$.  The compound Poisson process, $X(s)$, $s
\geq 0$ is defined by
\begin{equation}
 D(s) = \sum_{k=1}^{N(s)} \xi_i,
\end{equation}
where $N(s)$ is a Poisson process with rate $1$.  In this case, the L\'{e}vy
measure  $\Pi$ in (\ref{e:LKformula})  is given by the
distribution of the random variable $\xi_1$.  Therefore, the L\'{e}vy exponent for this  
process is
\begin{equation}
\phi(\lambda) = \int_0^\infty (1-e^{-\lambda x}) \nu(dx).
\end{equation}

Notice that $X(t) = 0$ for all $ t < \tau_1$, where $\tau_1$ has an
exponential distribution with mean $1$.  This implies that $E(0)$ will
also have an exponential distribution, and thus $E(0) \neq 0$ a.s.  The fact that $E(0) \neq 0$ is
characteristic of compound Poisson processes. Indeed, suppose that the inverse of a subordinator $D$ satisfies $E(0) > 0$ a.s., then in particular,
$\mathbb{E}E(t) \rightarrow c > 0$ at $t \rightarrow 0$ since $E$ has
cadlag paths.  Thus,  by the
Tauberian Theorem, 
\begin{equation}
\widetilde{U}(\lambda) = \frac{1}{\lambda \phi(\lambda)} \sim
{\cal{L}}[c](\lambda) = \frac{c}{\lambda}, \quad
\lambda \rightarrow \infty.
\end{equation}  
Thus, $\phi(\lambda)/c \rightarrow 1$ as $\lambda \rightarrow \infty$,
implying $\phi$ is bounded, which means $D$ is a Compound Poisson
process (see \cite{Bertoin:1996}, Corollary I.1.3).  

Now focus on the special case when $\xi_1$ has a Pareto distribution,
meaning $\xi_1$ has probability density
$\nu(x) = \alpha x^{-\alpha-1}$ for $x \geq 1$, and $\alpha
>0$ is fixed.  The L\'{e}vy exponent is given by
\begin{equation}\label{e:phiforpareto}
\phi(\lambda) = \alpha \int_1^\infty (1 - e^{-\lambda x}) x^{-\alpha-1} dx =
1 - \alpha \mbox{ \rm Ei}_{1+\alpha}(\lambda),
\end{equation}
where $\mbox{ \rm Ei}_{\alpha+1}(\lambda) = \int_1^\infty \frac{e^{-\lambda
    x}}{x^{\alpha+1}} dx$ is the Exponential integral.  Using a series
expansion\footnote{obtained using Mathematica} of $\mbox{ \rm Ei}_{\alpha+1}$, we have
\begin{equation}
\phi(\lambda) \sim \begin{cases} - \alpha \Gamma(-\alpha)
  \lambda^\alpha \quad 0 < \alpha < 1 \\
(1-\gamma_e-\log(\lambda))\lambda \quad \alpha =1 \\ 
\frac{\alpha}{\alpha -1} \lambda \quad \alpha > 1 \end{cases} \quad
\lambda \rightarrow 0,
\end{equation}
where $\gamma_e \approx 0.5772$ is the Euler constant.  Since
$\widetilde{U}(\lambda) = (\lambda \phi(\lambda))^{-1}$, the Tauberian
Theorem gives the large time behavior of the first passage time for
this Compound Poisson process
\begin{equation}\label{e:asypareto}
U(t) \sim \begin{cases}
 \displaystyle \frac{1}{-\alpha \Gamma(-\alpha) \Gamma(1+\alpha)} t^\alpha,\quad
  0<\alpha < 1 \\
\\
\displaystyle \frac{1}{1-\gamma_e + \log(t)} t \quad \alpha =1 \\
\\
\displaystyle \frac{\alpha-1}{\alpha} t ,\quad \alpha >1
\end{cases}, \quad t \rightarrow \infty.
\end{equation}

\begin{figure}[ht]
\centering
\includegraphics[width=0.8\textwidth]{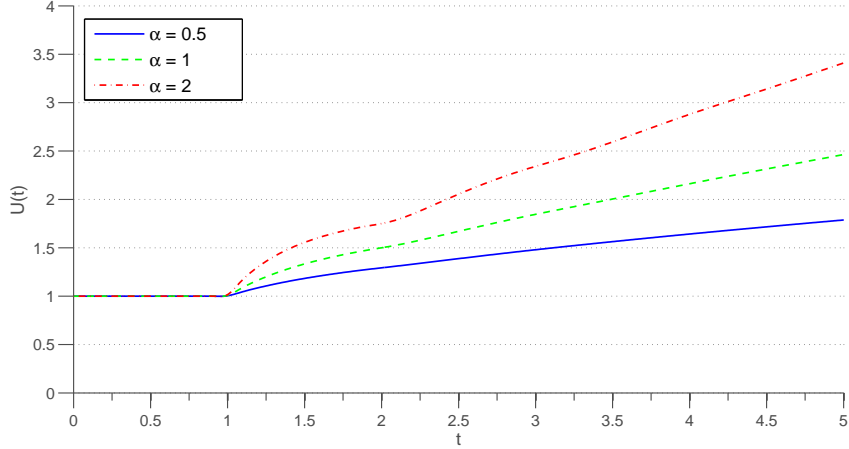}
\caption{  The mean
  first passage time of the compound Poisson process with Pareto jumps
  computed using the Post-Widder inversion method for $\alpha = 0.5,1,2$. }
\label{f:pareto}
\end{figure}

\begin{figure}[ht]
\centering
\includegraphics[width=0.8\textwidth]{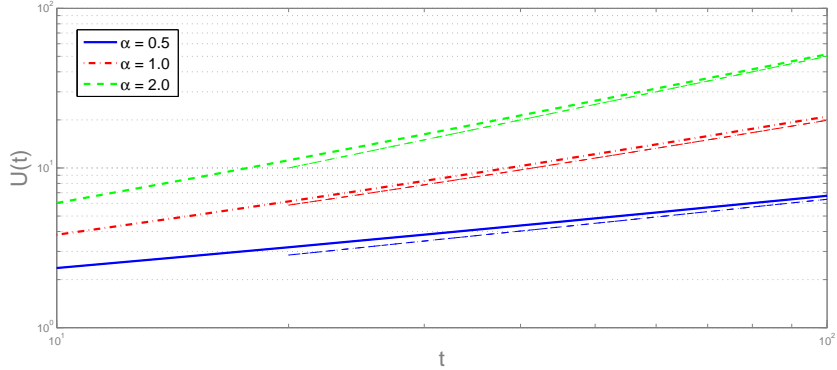}
\caption{  A log-log plot of the mean
  first passage time of the compound Poisson process with Pareto jumps
  computed using the Post-Widder inversion method for $\alpha =
  0.5,1,2$.  The dotted lines represent the asymptotic expressions
  given in (\ref{e:asypareto}).  We see good agreement for large times.   }
\label{f:paretolog}
\end{figure}

\begin{figure}[ht]
\centering
\includegraphics[width=0.8\textwidth]{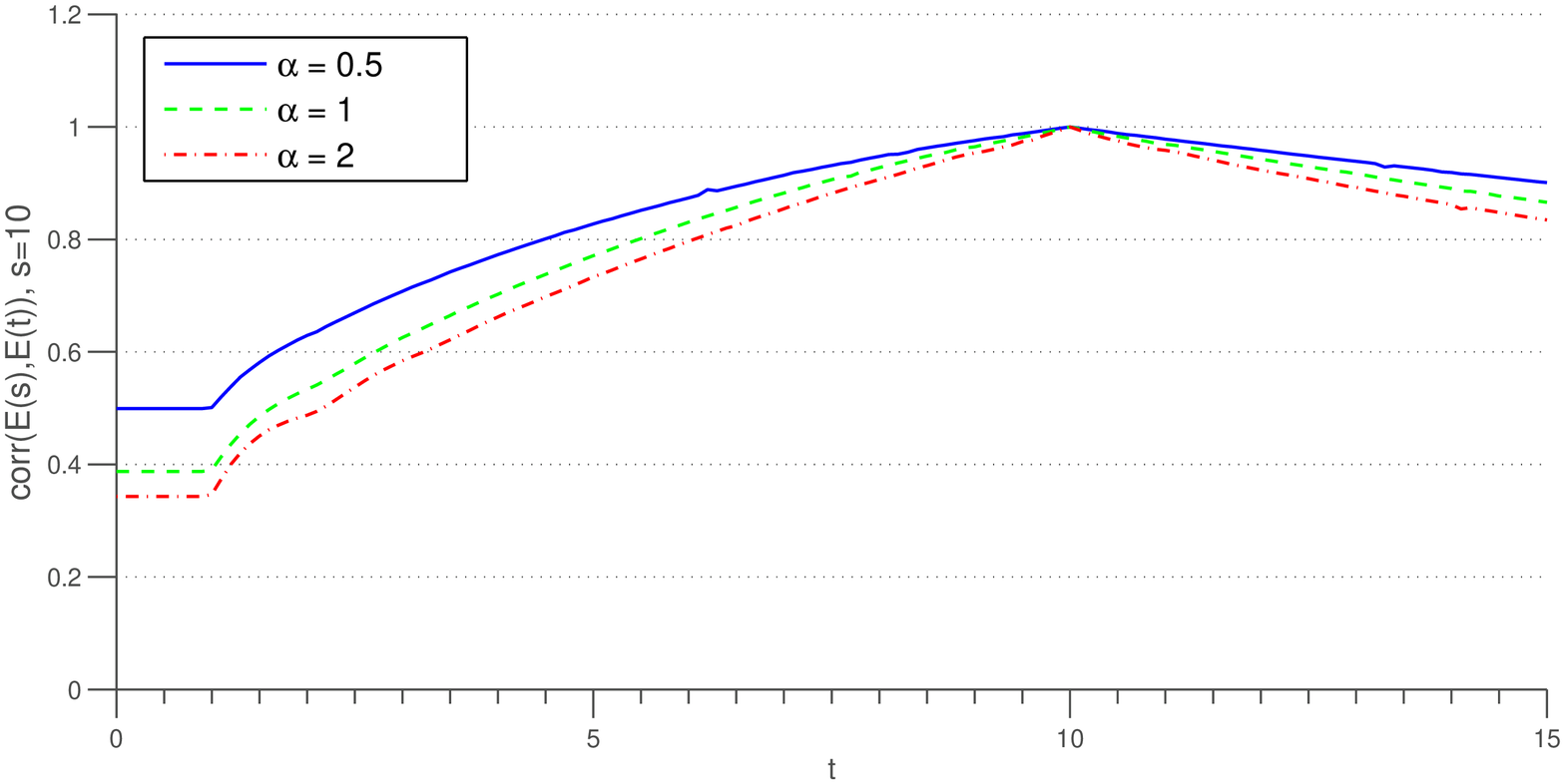}
\caption{  The correlation coefficient $\mbox{corr}(E(t),E(s))$ with
  $s=10$ fixed for
  the compound process with Pareto jumps for $\alpha = 0.5,1,2$. } 
\label{f:paretocorr}
\end{figure}

To calculate $U$, we found the Post-Widder method in
Section \ref{s:numerical} to be most useful\footnote{The only problem
 with this method occurs near $t=1$ where this technique is slow to converge.  The
  reason for this seems to be that $U$ is not differentiable here.  }.  In Figure \ref{f:pareto}
we plot $U$ for $\alpha = 1/2,1,2$, and in Figure \ref{f:paretolog} we
plot $U$ along with the asymptotic expressions given in
(\ref{e:asypareto}) with log-log scale.  The correlation
$\mathrm{corr}(E(t),E(s))$ is plotted in Figure \ref{f:paretocorr}.
Note that the renewal measure $dU$ will give a mass of weight $1$ to
$t=0$ since $U$ has a jump of size $1$ there.  For $t>0$, $dU$ is a
measure which has density $U'(t)$ which can be calculated as $U$ is
with the Post-Widder method as noted in a remark in Section \ref{s:postwidder}.

\subsection{``Mixture'' of $\alpha$-stable subordinators}\label{s:astable}

We consider here a continuous mixture of $\alpha$-stable subordinators
with $0 < \alpha < 1$.  Namely, the subordinator whose
L\'{e}vy exponent is given by
\begin{eqnarray}
\phi(\lambda) &=& \int_0^1 p(\beta) \lambda^\beta
d\beta  \label{e:mixed} \\
&=&  \int_0^\infty (1 - e^{-\lambda x}) g_p(x) dx,
\end{eqnarray}
where $p$ is a probability density on $(0,1)$ and the density $g_p$ of
the L\'{e}vy measure is
given by
\begin{equation}
g_p(x) = \int_0^1 \frac{x^{-\beta-1}}{\Gamma(-\beta)} p(\beta) d\beta
\end{equation}
In general, this
subordinator has no finite moments.  The $\alpha$-stable subordinator
corresponds to the choice $p(\beta;\alpha) = \delta(\beta-\alpha)$
in (\ref{e:mixed}).  Here we will consider two
extensions of this, namely when we choose a sum of two $\alpha$-stables,
$p(\beta;\alpha_1,\alpha_2) = C_1 \delta(\beta-\alpha_1) + C_2 \delta(\beta-\alpha_2)$,
with $\alpha_2 < \alpha_1$ and $C_1+C_2 = 1$, as well as what we call the ``uniform mix'',
which corresponds to the choice $p(\beta) = 1$ on $(0,1)$. 

\subsubsection{  Single $\alpha$-stable}  Many properties of the inverse of
an $\alpha$-stable subordinator with $0 < \alpha < 1$ are known (see for
instance, \cite{Baule:2005}), but we restate them here for
convenience.  Consider the $\alpha$-stable subordinator $\{D_\alpha(s),
\ s \geq 0\}$ with L\'{e}vy exponent $\phi(\lambda) = \lambda
^\alpha$.   Taking the inverse Laplace transform of (\ref{e:lt_mfpt}),
we have that the mean first passage
time for an inverse $\alpha$-stable subordinator is given by
\begin{equation}
U(t) = \frac{t^\alpha}{\Gamma(1+\alpha)},
\end{equation}
which implies that the density of the renewal measure is $U'(t) =
t^{\alpha - 1}/\Gamma(\alpha)$.  With this, the covariance can
be given in closed form.  Assume $s \leq t$, then 
\begin{eqnarray}
\mbox{\rm Cov}(E(s),E(t)) &=& \frac{1}{\Gamma(1+\alpha) \Gamma(\alpha)} \int_0^s
\left( (t-\tau)^\alpha + (s-\tau)^\alpha \right) \tau^{\alpha-1} d\tau
- \frac{(st)^\alpha}{\Gamma(1+\alpha)^2} \\
&=&  \frac{s^{2 \alpha}}{\Gamma(1+2\alpha)} + s^\alpha t^\alpha \left(
  \frac{1}{\Gamma(1+\alpha)}
  F(\alpha,-\alpha,\alpha+1;s/t) - \frac{1}{\Gamma(1+\alpha)^2}
\right).  \label{e:astablecov}
\end{eqnarray}
The above integral was computed in Mathematica and $F$ denotes the regularized confluent hypergeometric function.

\subsubsection{  Sum of two $\alpha$-stable subordinators}  
\begin{figure}[ht]
\centering 
\includegraphics[width=1.1\textwidth]{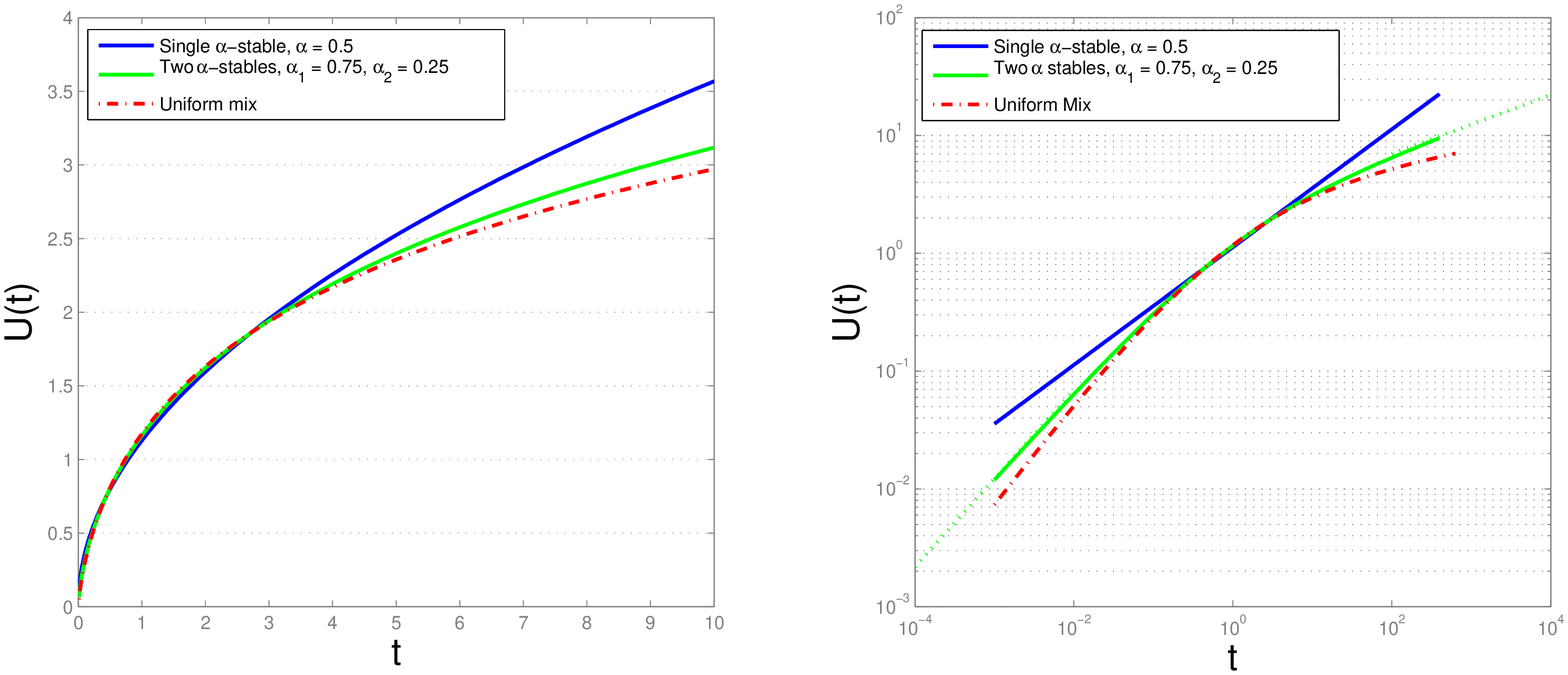}
\caption{ Plot of $U(t)$ calculated with the Post-Widder method for
  a single $\alpha$-stable, a sum of two $\alpha$ stables, and a
  uniform mixture.  On the right we used a
log-log plot and dotted lines to show asymptotic behavior given by
(\ref{e:asy}) for the sum of two $\alpha$-stable case.}
\label{f:sumalpha}
\end{figure}

Here we consider
the case when $p(\beta) = C_1\delta(\beta - \alpha_1) + C_2
\delta(\beta - \alpha_2)$ and the L\'{e}vy exponent is given by 
\begin{equation}\label{e:phiforsumas}
\phi(\lambda) =
C_1 \lambda^{\alpha_1} + C_2 \lambda^{\alpha_2},
\end{equation}
where $\alpha_2 < \alpha_1$ and $C_1+C_2 = 1$.
This corresponds to the subordinator given by $C_1^{1/\alpha_1} D_{\alpha_1}(s) +
C_2^{1/\alpha_2} D_{\alpha_2}(s)$.  From (\ref{e:lt_mfpt}), the Laplace
transform of the mean first passage time is given by
\begin{equation}\label{e:mfpt_two}
\widetilde{U}(\lambda) =  \frac{1}{ C_1
  \lambda^{\alpha_1+1} + C_2\lambda^{\alpha_2+1}}.
\end{equation} 
Using the Tauberian Theorem,  $U$ has the following
behavior in the limits $t \rightarrow 0$ and $t \rightarrow \infty$,
\begin{equation}\label{e:asy}
U(t) \sim \begin{cases} \displaystyle  \frac{t^{\alpha_1}}{C_1
    \Gamma(1+\alpha_1)}, \quad t \rightarrow 0 \\
\\
\displaystyle \frac{t^{\alpha_2}}{C_2
    \Gamma(1+\alpha_2)}, \quad t\rightarrow \infty. \end{cases}
\end{equation}
Thus, we see a cross-over in power-law behavior (which is displayed
with a log-log plot in Figure \ref{f:sumalpha}).  A closed form for the inverse Laplace transform of (\ref{e:mfpt_two})
could not be found, however $U$ can be calculated numerically using
the Post-Widder method (method 2 in Section \ref{s:numerical}).
Figure \ref{f:sumalpha} shows plots of this function for various
parameter values.  

  To obtain the correlation, we set $dU(t) = U'(t)
dt$, calculated $U'(t)$ using the Post-Widder method and then
evaluated the integrand in (\ref{e:inversecov1}) using numerical
integration.  The correlation is plotted in Figure
\ref{f:mixed_cor}.

\subsubsection{  Uniform mixture}

\begin{figure}[ht]
\centering
\includegraphics[width=0.8\textwidth]{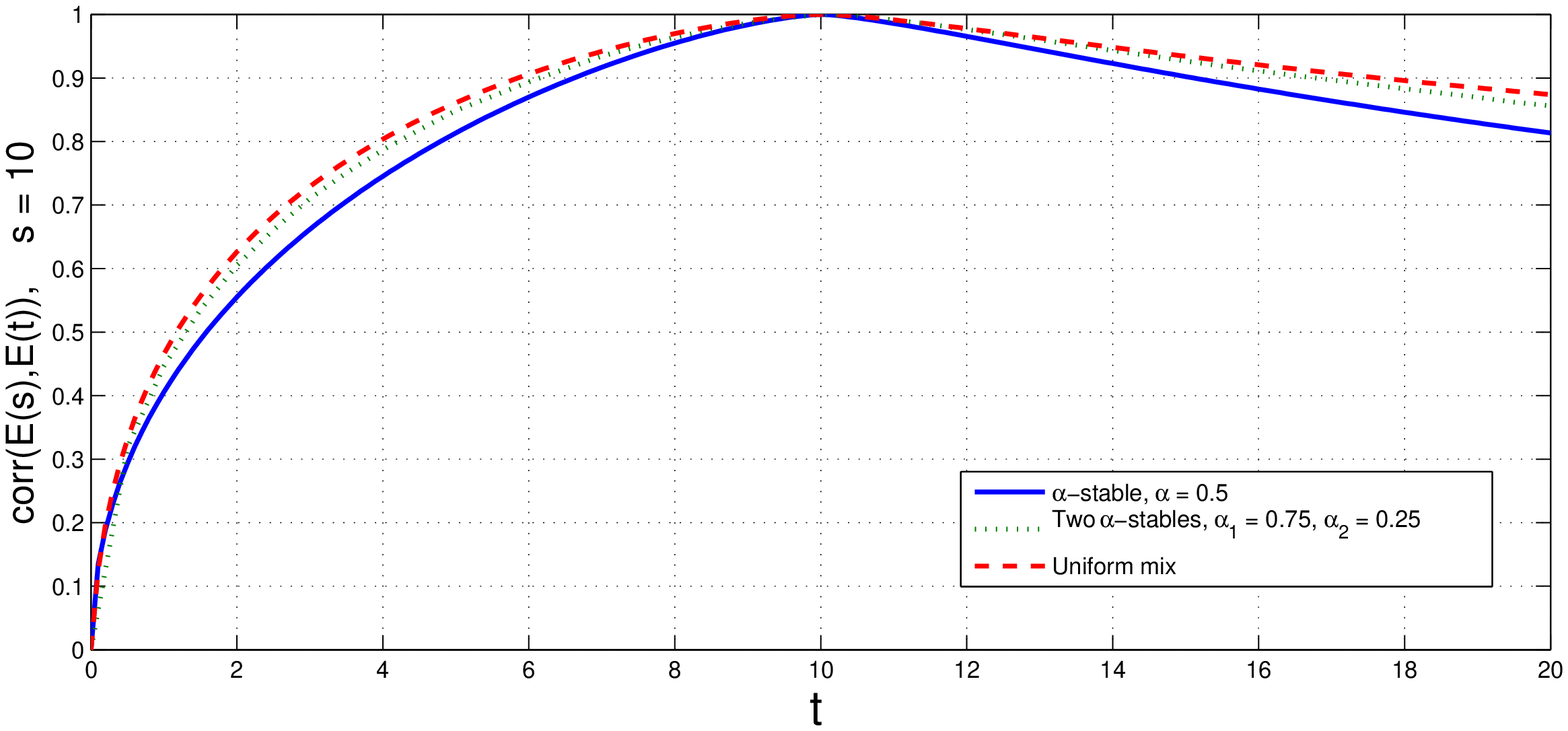}
\caption{A plot of the correlation coefficient with $s$ fixed at $10$
  and $t$ varying for the
  $\alpha$-stable process with $\alpha = 1/2$, the sum of two
  $\alpha$-stables with $\alpha_1 = .75$ and $\alpha_2 = .25$ and the
  uniform mix.  The correlation for the $\alpha$-stable case can be
  calculated exactly using (\ref{e:astablecov}).  In the other two
  cases, we numerically approximated the integral in (\ref{e:inversecov1}).}
\label{f:mixed_cor}
\end{figure}

Here we consider the case when $p(\beta) = 1$, for $\beta \in (0,1)$.
 The L\'{e}vy exponent is given by
\begin{equation}
\phi(\lambda) = \int_0^1 \lambda^\beta d\beta = \frac{\lambda -
  1}{\log(\lambda)}.
\end{equation}
The mean-first passage time for the inverse of this process has a closed form expression
given by (see \cite{Roberts:1966}, Eq. 4.1.9)
\begin{equation}
U(t) = {\cal L}^{-1} \frac{1}{\lambda \phi(\lambda)} =  {\cal{L}}^{-1} \frac{\log(\lambda)}{\lambda^2 - \lambda} = \gamma_e + e^t
\Gamma(0,t) +\log(t).
\end{equation}
Here $\Gamma(0,t)$ is the incomplete gamma function given by
$\Gamma(0,t) = \int_t^\infty e^{-z} z^{-1} dz$.  Since, for $t$ large, 
\begin{equation}
e^t\Gamma(0,t) = \int_t^\infty e^{-(z-t)}z^{-1} dz \leq \int_0^\infty
e^{-z} dz = 1,
\end{equation} 
we obtain the ``ultraslow'' growth:
\begin{equation} 
U(t) \sim \log(t), \quad t \rightarrow \infty.
\end{equation}
The density of the renewal measure also has a nice form:
\begin{equation}
U'(t) = e^t \Gamma(0,t) - e^t e^{-t} t^{-1} + t^{-1} =  e^t \Gamma(0,t).
\end{equation}
Using this, we can numerically calculate the covariance for the inverse
of the uniform mixed subordinator with (\ref{e:inversecov1}).  The correlation is plotted in
Figure \ref{f:mixed_cor}.

\subsection{Generalized inverse Gaussian L\'{e}vy processes}\label{s:gig}

The generalized inverse Gaussian (GIG) distribution \cite{barndorff:1977},~\cite{eberlein:2002} is a distribution characterized by three
parameters $\delta,\gamma,\kappa$ and has probability density function
given by
\begin{equation}\label{e:pforgig}
p_{GIG(\delta,\gamma,\kappa)}(x) =
\left(\frac{\gamma}{\delta}\right)^{\kappa} \frac{1}{2 K_\kappa(\delta
  \gamma)} x^{\kappa - 1} e^{-\frac{1}{2} ( \delta^2 x^{-1} + \gamma^2
  x)}, \quad x > 0.
\end{equation}
Here, $K_\kappa$ denotes the modified Bessel function of the third
kind (\cite{gray:1966}, section 3.2).  The parameters of this
distribution may take the following values:
\begin{eqnarray}
\delta \geq 0, \  \gamma > 0, & \mbox{if} & \ \kappa > 0, \\
\delta > 0, \  \gamma > 0,  & \mbox{if} &\ \kappa = 0, \\
\delta > 0, \  \gamma \geq 0,  & \mbox{if}& \ \kappa < 0.
\end{eqnarray} 
Important subclasses in this family are
\begin{itemize}
\item $\kappa > 0$, $\delta = 0$, $\gamma >0$ gives a Gamma distribution
  $\Gamma(\kappa,2/\gamma^2)$ with density
\begin{equation}\label{e:giggamma}
p_{GIG(0,\gamma,\kappa)}(x) = \frac{\gamma^{2 \kappa}}{2^{\kappa}\Gamma(\kappa)} x^{\kappa - 1}
e^{-\frac{\gamma^2}{2} x}, \quad x>0,\kappa>0,\gamma>0.
\end{equation}
\item $\kappa < 0$, $\delta > 0$, $\gamma = 0$ gives a reciprocal Gamma distribution
  $\mathrm{R}\Gamma(\kappa,\delta^2/2)$ with density
\begin{equation}\label{e:gigrgamma}
p_{GIG(\delta,0,\kappa)}(x) = \frac{\delta^{-2 \kappa}}{2^{-\kappa}
  \Gamma(-\kappa)} x^{\kappa-1} e^{-\frac{\delta^2}{2 x}}, \quad
x>0,\kappa<0,\delta>0.
\end{equation}
This distribution only has finite moments of order less than $|\kappa|$.

\item $\kappa = -\frac{1}{2}$, $\delta >0$,$\gamma \geq 0$ gives an
  Inverse Gaussian distribution $\mathrm{IG}(\delta,\gamma)$ with density
\begin{equation}
p_{GIG(\delta,\gamma,-1/2)}(x) = \frac{\delta e^{\gamma \delta}}{\sqrt{2 \pi}}
x^{-3/2}e^{ -\frac{1}{2}(\delta^2 x^{-1} + \gamma^2 x)},  \quad
x>0,\gamma \geq 0,\delta>0.
\end{equation}
The Inverse Gaussian distribution is the distribution of the first time a
Brownian motion with variance $\delta$ and drift $\gamma$ reaches the
level $1$ (\cite{Applebaum:2004}, Example 1.3.21).  All moments are finite if $\gamma > 0$, but if $\gamma=0$, it becomes
a totally right-skewed $\frac{1}{2}$-stable distribution which only has finite moments of order
less than $1/2$. 

\end{itemize}

The GIG distribution was shown to be infinitely divisible in
\cite{barndorffn:1977} and its L\'{e}vy-Khintchine representation was
derived in Section 5 of
\cite{eberlein:2002}.  The L\'{e}vy-Khintchine representation given in
\cite{eberlein:2002} is not in the form (\ref{e:LKformula}) and some
computation is needed to bring it into this form.  This is done in
appendix \ref{a:LKform}.  

The
Laplace transform of $p_{GIG(\delta,\gamma,\kappa)}$ for $\delta>0$
and $\gamma>0$ is given by
\begin{eqnarray}\label{e:gigmgf}
\widetilde{p}_{GIG(\delta,\gamma,\kappa)}(\lambda) &=& \left(
  \frac{\gamma^2}{\gamma^2 + 2 \lambda} \right)^{\kappa/2}
\frac{K_\kappa(\delta\sqrt{\gamma^2 + 2 \lambda})}{K_\kappa(\delta
  \gamma)} \\
&=& \exp(- \phi_{GIG(\delta,\gamma,\kappa)}(\lambda)),
\end{eqnarray}
where the L\'{e}vy exponent $\phi_{GIG}$ is given by
\begin{equation}
\phi_{GIG(\delta,\gamma,\kappa)}(\lambda) = \int_0^\infty (1 -
e^{-\lambda x}) g_{GIG(\delta,\gamma,\kappa)}(x) dx, \label{e:giglevy}
\end{equation}
Thus, the L\'{e}vy exponent is of the form (\ref{e:LKformula}) with
drift $\mu = 0$.  The corresponding L\'{e}vy measure (density) is (see
\cite{eberlein:2002} or Appendix \ref{a:LKform})
\begin{equation}\label{e:levymeasure}
g_{GIG(\delta,\gamma,\kappa)}(x) =  \frac{e^{-\frac{\gamma^2}{2} x}}{x}\left( \int_0^\infty
  \frac{e^{-x y}}{\pi^2 y (J^2_{|\kappa|}(\delta \sqrt{2 y})  +
    Y^2_{|\kappa|}(\delta \sqrt{2 y}))}dy + \max(0,\kappa) \right),
\quad x >0.
\end{equation} 
Here, $J_\nu$ and $Y_\nu$ denote the Bessel function of the first and
second kind, respectively, with index $\nu$ (\cite{gray:1966}, Chapter
2).  

By
letting $\delta \rightarrow 0$ with $\kappa >0$ and $\gamma
\rightarrow 0$ with $\kappa <0$ 
in (\ref{e:gigmgf}),  we obtain the
Laplace transform of the gamma distribution (\ref{e:giggamma}) and the reciprocal gamma
distribution (\ref{e:gigrgamma}) respectively.  Doing so gives (see \cite{eberlein:2002})
\begin{eqnarray}
\widetilde{p}_{GIG(0,\gamma,\kappa)}(\lambda) &=& \left( 1 + \frac{2
    \lambda}{\gamma^2} \right)^{-\kappa} \quad \kappa,\gamma > 0,\label{e:gammamgf} \\ 
\widetilde{p}_{GIG(\delta,0,\kappa)}(\lambda) &=&
\frac{2\left(\frac{\delta^2}{2}\right)^{-\kappa/2}
  \lambda^{-\kappa/2}}{\Gamma(-\kappa)} K_{-\kappa}(\delta \sqrt{2
  \lambda}), \quad \kappa<0,\delta >0. \label{e:rgammamgf}
\end{eqnarray}
The corresponding L\'{e}vy measures are obtained similarly using
(\ref{e:levymeasure}), (recall that $Y_{\nu}(z) \rightarrow -\infty$
as $z \rightarrow 0$):
\begin{eqnarray}
g_{GIG(0,\gamma,\kappa)}(x) &=& \frac{\kappa}{x} e^{-\frac{\gamma^2}{2}
  x}  \label{e:gigmeas1} \\
g_{GIG(\delta,0,\kappa)}(x) &=& \frac{1}{x}  \int_0^\infty
  \frac{e^{-x y}}{\pi^2 y (J^2_{|\kappa|}(\delta \sqrt{2 y})  +
    Y^2_{|\kappa|}(\delta \sqrt{2 y}))}dy.  \label{e:gigmeas2} 
\end{eqnarray}

The mean of a GIG distribution, when it exists, can be calculated
from these Laplace transforms.  If $X \sim GIG(\delta,\gamma,\kappa)$,
then there are three cases for $\mathbb{E}X$:
\begin{eqnarray}
\mathbb{E}X &=& \frac{\delta K_{1+\kappa}(\gamma \delta)}{\gamma
  K_\kappa(\gamma \delta)}, \quad \kappa \in
\mathbb{R},\gamma>0,\delta>0 \label{e:gigmean1} \\
\mathbb{E}X &=& \frac{2 \kappa}{\gamma^2},\quad
\delta=0,\kappa>0,\gamma>0 \label{e:gigmean2} \\
\mathbb{E}X &=& \begin{cases} \infty, \quad \gamma = 0, \quad -1 \leq \kappa < 0, \delta >0  \\
 \frac{1}{2(-\kappa - 1) }
  \delta^{2}, \quad \quad  \quad \gamma = 0, \kappa < -1, \delta >0 \end{cases}
\end{eqnarray}

We talked so far about the GIG distribution.  We now consider the
corresponding L\'{e}vy process, namely, the subordinator $\{D_{GIG(\delta,\gamma,\kappa)}(s), \
s \geq 0 \}$ with Laplace transform
\begin{equation}
\mathbb{E}\exp(-\lambda D_{GIG(\delta,\gamma,\kappa)}(s)) =
\exp(-s\phi_{GIG(\delta,\gamma,\kappa)}(\lambda)).
\end{equation}
From  (\ref{e:giglevy}),
the drift of this process is $0$ and its L\'{e}vy measure is given by
$\Pi(dx) = g_{GIG}(x) dx$.  Notice that this L\'{e}vy process is
indeed a subordinator since its L\'{e}vy measure is concentrated on
the positive axis and from (\ref{e:gigmeas1}) and (\ref{e:gigmeas2}),
it follows that $\int_0^\infty (1 \wedge x) \Pi(dx) < \infty$.  Let $E_{GIG(\delta,\gamma,\kappa)}$ be the
inverse of $D_{GIG(\delta,\gamma,\kappa)}$ and let
$U_{GIG(\delta,\gamma,\kappa)}(t) = \mathbb{E}E_{GIG(\delta,\gamma,\kappa)}(t)$.   The gamma process
($\delta = 0$) and the inverse Gaussian process ($\kappa = -1/2$) are
treated in \cite{Nord:2005}, where closed form expressions are given
for the renewal measure.  In general, we cannot write $U_{GIG}(t)$ in
closed form, however, we can obtain
expressions for large and small time behavior using the
Tauberian theorem. 

\medskip

\noindent \textbf{Asymptotics of $U_{GIG(\delta,\gamma,\kappa)}$, $t
  \rightarrow 0$}:  One has $K_\kappa(x) \sim \sqrt{\frac{\pi}{2
    x}} e^{-x}$ as $x \rightarrow
\infty$ ( \cite{gray:1966}, section 5.4).  We first consider the case
when $\delta > 0$.  When $\gamma > 0$,  (\ref{e:gigmgf}) implies that as $\lambda \rightarrow \infty$,
\begin{eqnarray}
\phi_{GIG(\delta,\gamma,\kappa)}(\lambda) &=& - \log\left(\left(
  \frac{\gamma^2}{\gamma^2 + 2 \lambda} \right)^{\kappa/2}
\frac{K_\kappa(\delta\sqrt{\gamma^2 + 2 \lambda})}{K_\kappa(\delta
  \gamma)} \right) \\
&=& \frac{\kappa}{2} \log\left( 1 + \frac{2 \lambda}{\gamma^2}\right)
- \log(K_k(\delta \sqrt{\gamma^2 + 2 \lambda})) + \log(K_\kappa(\delta
  \gamma)) \\
&=& \delta \sqrt{ \gamma^2 + 2 \lambda} + O(\log(\lambda)) \\
& \sim &  \delta \sqrt{2 \lambda}, \quad \lambda \rightarrow \infty.
\end{eqnarray}
The case $\gamma = 0$ gives the same answer (by using (\ref{e:rgammamgf})).  If $\delta = 0$, we instead use (\ref{e:gammamgf}) and get
\begin{eqnarray}
\phi_{GIG(0,\gamma,\kappa)}(\lambda) &=& \kappa \log\left(  1 + \frac{2
    \lambda}{\gamma^2} \right) \\
&\sim& \kappa \log( \lambda), \quad \lambda \rightarrow \infty.
\end{eqnarray}

Since $\widetilde{U}_{GIG(\delta,\gamma,\kappa)}(\lambda) = (\lambda
\phi_{GIG(0,\gamma,\kappa)}(\lambda))^{-1}$, we get from the Tauberian theorem
\begin{equation}\label{e:gigasyzero}
U_{GIG(\delta,\gamma,\kappa)}(t) \sim \begin{cases}
 \displaystyle  \frac{-1}{\kappa \log(t)}, \quad \delta = 0, \\
\\
\displaystyle \sqrt{\frac{2}{\pi \delta^2}} \sqrt{t}, \quad \delta>0 \end{cases} \quad t
\rightarrow 0.
\end{equation}

\medskip

 \textbf{Asymptotics of $U_{GIG(\delta,\gamma,\kappa)}$, $t
  \rightarrow \infty$}:  If $\gamma > 0$, the renewal theorem implies
that $U_{GIG(\delta,\gamma,\kappa)}(t) \sim t/\mathbb{E}D(1)$, where
$\mathbb{E}D(1)$ is finite and is given by (\ref{e:gigmean1}) or (\ref{e:gigmean2}). 

 The case $\gamma = 0$ requires more work. We shall use the following series expansion of
$K_{\kappa}(x)$ as $x \rightarrow 0$ with $\kappa>0$ (see \cite{Mathai:1993}, page 121):
\begin{equation}\label{e:besselasy}
K_{\kappa}(x) = \begin{cases}
 \displaystyle \frac{\Gamma(\kappa)}{2^{1-\kappa}} x^{-\kappa} \left( 1 +
    \frac{\Gamma(-\kappa)}{4^\kappa \Gamma(\kappa)} x^{2 \kappa}
  \right) + o(x^\kappa), \quad 0 < \kappa < 1 \\
\\
\displaystyle x^{-1}\left( 1 + \frac{1}{4}(2 \gamma_e - 1 + 2 \log(x) -
  \log(4))x^2\right)  + o(x), \quad \kappa = 1 \\
\\
\displaystyle \frac{\Gamma(\kappa)}{2^{1-\kappa}} x^{-\kappa} \left( 1
  + \frac{1}{4 (1 - \kappa)}x^2\right) + o(x^{2-\kappa}), \quad \kappa > 1
\end{cases}.
\end{equation}

 Now, (\ref{e:rgammamgf}) and (\ref{e:besselasy})
imply that for $0 < -\kappa < 1$,
\begin{eqnarray}
\phi_{GIG(\delta,0,\kappa)}(\lambda) &=& - \log\left(  \frac{2\left(\frac{\delta^2}{2}\right)^{-\kappa/2}
  \lambda^{-\kappa/2}}{\Gamma(-\kappa)} K_{-\kappa}(\delta \sqrt{2
  \lambda})  \right) \\
&=& - \log \left( \frac{2^{1+\kappa/2} \delta^{-\kappa}
    \lambda^{-\kappa/2}}{\Gamma(-\kappa)}\right) - \log \left(
  K_{-\kappa}(\delta \sqrt{2 \lambda} ) \right) \\
&=& -\log\left( \frac{2^{1+\kappa/2} \delta^{-\kappa}
    \lambda^{-\kappa/2}}{\Gamma(-\kappa)}\right) -
\log\left(\frac{\Gamma(-\kappa)}{2^{1+\kappa}} \delta^\kappa
  2^{\kappa/2} \lambda^{\kappa/2}\right) - \log\left( 1 +
  \frac{\Gamma(\kappa)}{2^{-\kappa} \Gamma(-\kappa)} \delta^{-2 \kappa}
   \lambda^{-\kappa} + o(\lambda^{-\kappa}) \right) \nonumber \\
&=& - \log\left( 1 +
  \frac{\Gamma(\kappa)}{2^{-\kappa} \Gamma(-\kappa)} \delta^{-2 \kappa}
   \lambda^{-\kappa} + o(\lambda^{-\kappa}) \right) \nonumber \\
&\sim& -\frac{\Gamma(\kappa)}{2^{-\kappa} \delta^{2 \kappa}
  \Gamma(-\kappa)} \lambda^{-\kappa}, \quad \lambda \rightarrow 0.
\end{eqnarray}
Similar calculations give for $\kappa = -1$,
\begin{eqnarray}
\phi_{GIG(\delta,0,\kappa)}(\lambda) &=& -\log \left(  \delta \sqrt{ 2
    \lambda} K_1(\delta \sqrt{2 \lambda}) \right) \\
&=& -\log\left( \delta \sqrt{ 2 \lambda} \right) - \log( K_1(\delta \sqrt{2
    \lambda}) ) \\ 
&=& -\log \left( \delta \sqrt{2 \lambda} \right) - \log\left( (\delta
  \sqrt{2 \lambda})^{-1} \right) - \log\left( 1 + \frac{1}{4} ( 2 \gamma_e - 1 +
 2 \log(\delta \sqrt{2 \lambda}) - \log(4)) ( 2 \delta^2 \lambda) +
 o(x^2) \right) \nonumber \\
&\sim& -\frac{\delta^2}{2}(2 \gamma_e - 1 +
 \log(2 \delta^2) -\log(4) + \log(\lambda)) \lambda,  \quad \lambda
 \rightarrow 0,
\end{eqnarray}
and for $-\kappa > 1$, 
\begin{eqnarray}
\phi_{GIG(\delta,0,\kappa)}(\lambda)  &=& - \log\left(  \frac{2\left(\frac{\delta^2}{2}\right)^{-\kappa/2}
  \lambda^{-\kappa/2}}{\Gamma(-\kappa)} K_{-\kappa}(\delta \sqrt{2
  \lambda})  \right) \\
&=& - \log \left( \frac{2^{1+\kappa/2} \delta^{-\kappa}
    \lambda^{-\kappa/2}}{\Gamma(-\kappa)}\right) - \log \left(
  K_{-\kappa}(\delta \sqrt{2 \lambda} ) \right) \\ 
&=&  -\log \left( \frac{2^{1+\kappa/2} \delta^{-\kappa}
    \lambda^{-\kappa/2}}{\Gamma(-\kappa)}\right) -
\log\left(\frac{\Gamma(-\kappa)}{2^{1+\kappa}} \delta^\kappa
  2^{\kappa/2} \lambda^{\kappa/2}\right) - \log \left( 1 +
  \frac{1}{4(1+\kappa)} (2 \delta^2 \lambda)
   + o(\lambda) \right) \nonumber \\
&\sim& - \frac{\delta^2}{2
  (1+\kappa)} \lambda, \quad \lambda \rightarrow 0.
\end{eqnarray}

Thus, using the Tauberian theorem, we have the follow large time
behavior of the mean first passage time of the GIG process
\begin{equation}\label{e:gigasyinf}
U_{GIG(\delta,\gamma,\kappa)}(t) \sim \begin{cases} \displaystyle  \frac{\gamma
  K_\kappa(\gamma \delta)}{\delta K_{1+\kappa}(\gamma \delta)} t,
\quad \kappa \in \mathbb{R},\gamma \geq 0,\delta >0, \\
\\
\displaystyle -\frac{2^{-\kappa} \delta^{2 \kappa}
  \Gamma(-\kappa)}{\Gamma(\kappa) \Gamma(1-\kappa)} t^{-\kappa} ,\quad
\gamma = 0, -1 < \kappa < 0, \delta > 0 \\
\\
\displaystyle \frac{\delta^2/2}{2 \gamma_e - 1 +
 \log(2 \delta^2) -\log(4) + \log(t)} t, \quad \gamma = 0, \kappa =
1,\delta >0 \\
\\
\displaystyle \frac{2 (- \kappa-1)}{\delta^2} t, \quad \gamma = 0,\kappa < -1,
\delta > 0  \end{cases}, \quad t \rightarrow \infty
\end{equation}

\begin{figure}[ht]
\centering
\includegraphics[width=0.8\textwidth]{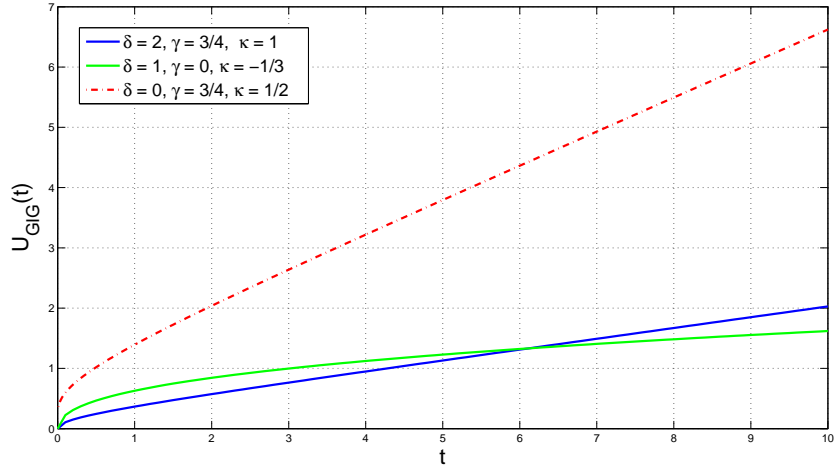}
\caption{Plots of $U_{GIG}$ for various values of parameters
  $\delta,\gamma$ and $\kappa$.  Each of these plots were generated
  using the Post-Widder method with $\epsilon = 10^{-5}$.}
\label{f:gig}
\end{figure}

\begin{figure}[ht]
\centering
\includegraphics[width=1.0\textwidth]{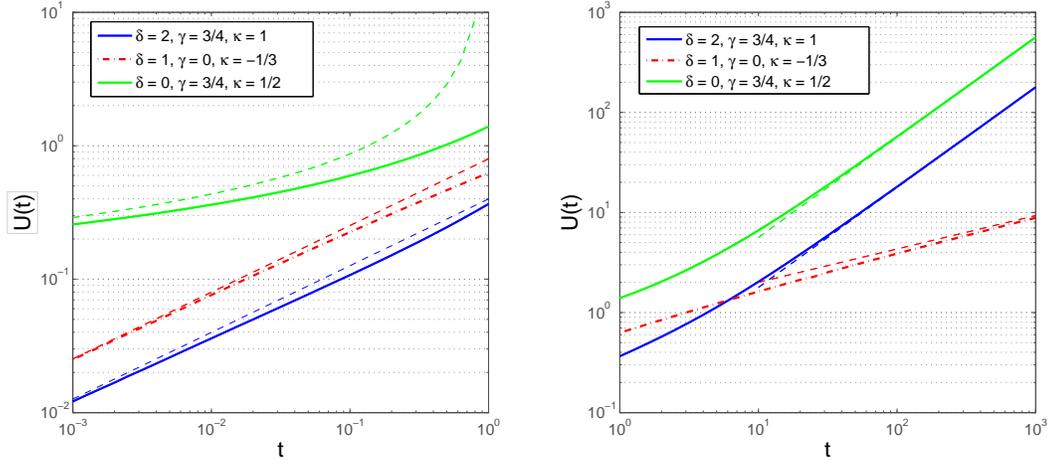}
\caption{Plots of $U_{GIG}$ on log-log scales.  To dotted curves
  correspond to the asymptotic expressions given by
  (\ref{e:gigasyzero}) (left), namely $ t \rightarrow 0$, and
  (\ref{e:gigasyinf}) (right), namely, $t \rightarrow \infty$.  We see
in these limits, the approximations are very good.}
\label{f:giglog}
\end{figure}

\begin{figure}[ht]
\centering
\includegraphics[width=0.8\textwidth]{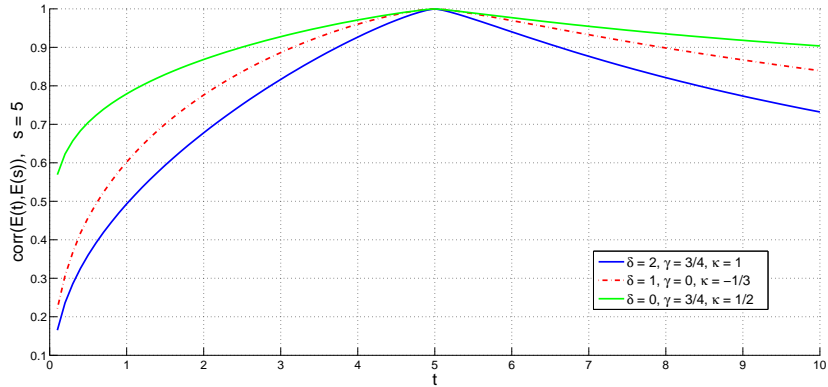}
\caption{Plots of the correlation coefficient $\rho =
  \mbox{corr}(E(t),E(s))$ for the generalized inverse Gaussian
  process with various parameters values.  Here, $s=5$ is fixed and
  $t$ varies}
\label{f:gigcorr}
\end{figure}

While (\ref{e:gigasyzero}) and (\ref{e:gigasyinf}) give the asymptotic of
$U$, closed form expressions for $U_{GIG}$ are not known.  Using our
methodology, we can compute $U_{GIG}$ numerically.  In Figure
\ref{f:gig}, $U_{GIG}$ is plotted for three sets of parameter values.
Figure \ref{f:giglog} shows a log-log plot which includes the asymptotic curves
(\ref{e:gigasyzero}) and (\ref{e:gigasyinf}) .  Since
this subordinator is strictly increasing, the renewal measure is given
by $dU(t) = U'(t)dt$, and $U'(t)$ can also be computed using the
Post-Widder approach.  The correlation $\mathrm{corr}(E(t),E(s))$ is
can also be calculated as it was in the other examples and is plotted in Figure \ref{f:gigcorr} with $s=10$ fixed.

%\textbf{$\alpha$-stable Process, $0 < \alpha < 1$}.  This example has attracted% much attention in the study of scaling limits of Continuous-time random walks,% where a Markovian process is subordinated to an inverse $\alpha$-stable proc%ess.  The use of $\alpha$-stable subordinators provides a model of anomalous d%iffusion, hence is useful in many physical applications.  Much of the results% of the previous sections have already been derived for this particular subordi%nator, see for instance.  For completeness we will present these results here %as simple applications of the results in this paper.    

\newpage

\section{Conclusion}

In this paper, we developed two numerical methods for calculating
the function $U(t) = \mathbb{E}E(t)$, where the process $\{ E(t)$, $t\geq
0 \}$ is the first hitting time of a  L\'{e}vy
subordinator $\{ D(s), s \geq 0\}$ with L\'{e}vy exponent $\phi$ given
by (\ref{e:LKformula}).  The
function $U$ has been shown to characterize all finite-dimensional
distributions of the process $E$ and is useful for calculating moments of
$E$, for example, the covariance (see equation (\ref{e:inversecov1})).  

The Laplace transform of $U$ has a simple expression in terms of
$\phi$, namely,  $\widetilde{U}(\lambda) = (\lambda
\phi(\lambda))^{-1}$.  Thus, calculating $U$ involves computing the
inverse Laplace transform of this function.  The
first method described in Section \ref{s:integration} computes the
inverse of this Laplace transform by approximating the Bromwich
integral given by (\ref{e:bromwich}).  This integral can be computed by
rewriting the integrand in terms of the real and imaginary parts of
$\phi$.  We give explicit expressions for $\mathrm{Re}( \phi)$ and
$\mathrm{Im}(\phi)$ in terms of the drift $\mu$ and the L\'{e}vy
measure $\Pi$ of the subordinator $D$.  

The second method described in Section \ref{s:postwidder} computes the
inverse Laplace transform of $\widetilde{U}$ using the Post-Widder inversion
formula given in (\ref{e:pwinversion}).  This formula is usually
difficult to use because it requires evaluating derivatives of high
orders.  However, using our methods, $\widetilde{U}^{(k)}$ can be calculated
in a reasonable fashion.  As with the integration case, all terms in
this approximation are given in terms of only the drift $\mu$ and
L\'{e}vy measure $\Pi$ corresponding to the subordinator $D$.   We tested both of these methods in cases
where $U$ can be calculated exactly and obtained accurate approximations.

As an application of our methods, we considered three families of 
L\'{e}vy subordinators $D$ and calculated the mean and correlation of
their respective inverse subordinators $E$.  The three families
considered were (i) Poisson and Compound Poisson processes (ii)
Continuous mixtures of $\alpha$-stable subordinators and (iii)
Generalized inverse Gaussian L\'{e}vy processes.  In each example, either
the integration or Post-Widder method was useful for calculating
the mean $U(t) = \mathbb{E}E(t)$.  Once we computed $U$,
the correlation function of $E$ can be computed by numerically
approximating the integral in (\ref{e:inversecov1}).   Along with these numerical
approximations, we also gave in each case, asymptotic expressions for $U$.

\section{Guide to Software}\label{s:guide}

We have developed a MATLAB software package which computes $U$ for the
examples above, as well as a program for a user defined example, which
is available from the authors.  The package includes 4 programs which calculate
$U(t)$:
\begin{itemize}

\item The Poisson process(\texttt{invert\_poisson.m})

\item Compound Poisson process with pareto
jumps (\texttt{invert\_pareto.m}), 

\item Sum of two $\alpha$-stable processes(\texttt{invert\_sumas.m})  

\item The generalized inverse
Gaussian L\'{e}vy process(\texttt{invert\_gig.m}).  

\end{itemize}

The density of the renewal measure $U'(t)$
is calculated in all cases but the Poisson process.  To use these
functions, invoke MATLAB and add the \texttt{inversesub}  directory to MATLAB's working path by typing 

\medskip

\noindent $>>$ \texttt{addpath('/yourpath/inversesub')}

\medskip

\noindent Here,  "yourpath" is the path in which the directory
\texttt{inversesub} is located  (for example, in Windows, this might
look something like ``\texttt{C:/myhomedir/inversesub}''). 

 Table \ref{t:guide} defines the
required inputs for each function including the method used (numerical
integral or Post-Widder), the required parameters, and the outputs one
obtains.   In each case, the default tolerance is set to $\epsilon = 10^{-6}$.
To change this, simply add an optional argument to each function which
gives the desired tolerance.  For example, typing

\medskip

\noindent $>>$  \texttt{invert\_poisson(1.2,0,1)}

\medskip

\noindent gives the mean first-hitting time $U(t) = \mathbb{E}E(t)$
for a Poisson process at time $t = 1.2$ with
drift $\mu =0$ a rate $r = 1$, and a tolerance of $10^{-6}$.  
Alternatively, one can type

\medskip

\noindent $>>$  \texttt{invert\_poisson(1.2,0,1,.001)}

\medskip  

\noindent to compute $U$ instead with a tolerance of $10^{-3}$.  
If the requested tolerance cannot be met, the program will return a
message saying so as well as a crude estimate of the error.

The programs using the Post-Widder method also returns an optional estimate for the
derivative of $U$, $U'(t)$.  For example, typing
\medskip

\noindent $>>$ \texttt{ [U DU] = invert\_gig(1,1,0,-1/2) }

\medskip   

\noindent assigns the value $U(1)$ to $U$ and $U'(1)$ to $DU$ where
$U$ corresponds to the inverse of the
reciprocal gamma L\'{e}vy process.  Here $t=1$ and $1,0,-1$  are
parameter values (see Table \ref{t:guide}).

Each program also accepts vector inputs for $t$.  For instance, 
\medskip

\noindent $>>$ \texttt{ [U DU] = invert\_gig([1 2 3],1,0,-1/2)}

\medskip    
\noindent assigns the vector $[U(1),\ U(2),\ U(3)]$ to the variable U and
the vector $[U'(1),\ U'(2),\ U'(3)]$ to the variable DU.  

The program ``invert\_empty.m'' contains all the code of the previous
examples with the piece which computes $\phi^{(n)}$ missing.  To use
this program, add a function call to line 30 of the code which
computes the derivatives of $\phi$ for your case.

\begin{table}[ht]
\centering 
\begin{tabular}{|l|l|l|l|l|}
\hline 
Function Name & Method Used & Number of Inputs & Meaning \& order of inputs &  Outputs \\ \hline
\multirow{3}{*} & & & $t$ :
time $t$ in $U(t)$ & \\
{\texttt{invert\_poisson}} &  {Numerical Integration}  &  {3} & mu : Drift $\mu$ in
 (\ref{e:phiforpois})  & {$U(t)$}  \\
 & & & r : Rate $r$ in (\ref{e:phiforpois})  &  \\ \hline
\multirow{2}{*}{\texttt{invert\_pareto}}  & {Post-Widder} & {2}  & $t$ :
time $t$ in $U(t)$ & $U(t),U'(t)$ \\
 & & & a : Exponent $\alpha$ in (\ref{e:phiforpareto})  &  \\ \hline
\multirow{5}{*}{\texttt{invert\_sumas}} & & & $t$ :
time $t$ in $U(t)$ & \\
 & & & a1 : $\alpha_1$ in (\ref{e:phiforsumas}) & \\ 
 &  {Post-Widder}  &  {5} & a2 :  $\alpha_2$ in (\ref{e:phiforsumas}) & {$U(t),U'(t)$}  \\
 & & & c1 : $C_1$ in  (\ref{e:phiforsumas})  &  \\
& & & c2 : $C_2$ in  (\ref{e:phiforsumas})  &  \\\hline
\multirow{4}{*} & & & $t$ :
time $t$ in $U(t)$ & \\
{\texttt{invert\_gig}} &  {Post-Widder}  & {4} & delta : $\delta$ in
(\ref{e:pforgig}) & {$U(t),U'(t)$} \\ 
 &  &   & gamma :  $\gamma$ in (\ref{e:pforgig}) &   \\
 & & & kappa : $\kappa$ in  (\ref{e:pforgig})  &  \\\hline
\end{tabular}
\caption{Information about the 4 functions included in the
  software package.  The parameters in the fourth column should be
  entered in the order of top to bottom, for example, for the sum of
  two $\alpha$-stable case, one would type
 \texttt{ invert\_sumas(t,a1,a2,c1,c2)} . } 
\label{t:guide}
\end{table}

\newpage

\appendix{

\section{Appendix:  L\'{e}vy-Khintchine form for GIG distributions} \label{a:LKform}

Here, we show that the L\'{e}vy exponent corresponding to a
GIG distribution can be written in the form
(\ref{e:LKformula}) with drift $\mu = 0$.  This is done in the following proposition. 

\begin{prop}
Let $p_{GIG(\delta,\gamma,\kappa)}$ denote the probability density of
the $GIG(\delta,\gamma,\kappa) $ distribution.  Then the Laplace
transform of $p_{GIG}$ is given by
\begin{equation}
\widetilde{p}_{GIG(\delta,\gamma,\kappa)}(\lambda) = \exp(-
\phi_{GIG(\delta,\gamma,\kappa)}(\lambda)), \label{e:LKF1}
\end{equation}
where the L\'{e}vy exponent $\phi_{GIG}$ is 
\begin{equation}\label{e:subform}
\phi_{GIG(\delta,\gamma,\kappa)}(\lambda) = \int_0^\infty (1 -
e^{-\lambda x}) g_{GIG(\delta,\gamma,\kappa)}(x) dx, 
\end{equation}
with L\'{e}vy measure (density)
\begin{equation}\label{e:lmeasure}
g_{GIG(\delta,\gamma,\kappa)}(x) = \begin{cases} \displaystyle \frac{e^{-\frac{\gamma^2}{2} x}}{x}\left( \int_0^\infty
  \frac{e^{-x y}}{\pi^2 y (J^2_{|\kappa|}(\delta \sqrt{2 y})  +
    Y^2_{|\kappa|}(\delta \sqrt{2 y}))}dy + \max(0,\kappa) \right),
\quad x >0, \delta > 0 \\
\displaystyle \frac{\kappa}{x} e^{-\frac{\gamma^2}{2} x} \quad x
\geq 0, \delta = 0 \end{cases}.
\end{equation} 
\end{prop}

\begin{pf}
 From \cite{eberlein:2002} Section 5.2, $\widetilde{p}_{GIG}$ is given
 in the form (\ref{e:LKF1}), but
 $\phi_{GIG}$ is expressed in following alternative L\'{e}vy
 representation which depends on the parameters $\delta, \gamma,\kappa$,
\begin{equation}\label{e:oldlf}
- \phi_{GIG(\delta,\gamma,\kappa)} (\lambda)  = \begin{cases} \displaystyle i \lambda
  \frac{\delta K_{1+\kappa}( \gamma \delta) }{\gamma K_{\kappa}(\gamma
    \delta)} + \int_0^\infty (e^{i \lambda x} - 1 - i \lambda x)
  g_{GIG(\delta,\gamma,\kappa)}(x) dx, \quad \delta, \gamma > 0, \kappa
  \in \mathbb{R} \\
\\
\displaystyle i \lambda \frac{2 \kappa}{\gamma^2} + \int_0^\infty (e^{i \lambda x} -
1 -i \lambda x) g_{GIG(0,\gamma,\kappa)}(x) dx,\quad \delta = 0, \kappa,\gamma>0 \\
\\
 \displaystyle i \lambda \delta^2 \int_0^\infty \frac{1 - e^{-x}}{x} g_{|\kappa|}(2
\delta^2 x) dx \\
\quad \displaystyle + \int_0^\infty (e^{i \lambda x} - 1 - i \lambda x
\bone_{[0,1]}(x) )g_{GIG(\delta,0,\kappa)}(x) dx,  \quad
\gamma=0, \kappa < 0, \delta>0
\end{cases}
\end{equation}
where, for $\nu > 0$,
\begin{equation}\label{e:defofg}
g_\nu(x) = \frac{2}{\pi^2 x [J_\nu^2 ( \sqrt{x}) + Y_\nu^2(\sqrt{x})
  ]}, \quad x > 0.
\end{equation}
Thus, to show (\ref{e:subform}), we much check that in each of the three cases in (\ref{e:oldlf}), the drift term
cancels with the ``$i \lambda x$'' term in the integrand.

\noindent \textbf{Case 1:} $\delta,\gamma>0, \kappa \in \mathbb{R}$.  

We must show
\begin{equation}\label{e:showthis1}
 \frac{\delta K_{1+\kappa}( \gamma \delta) }{\gamma K_{\kappa}(\gamma
    \delta)} = \int_0^\infty x g_{GIG(\delta,\gamma,\kappa)} dx.
\end{equation}
Using (\ref{e:lmeasure}) and performing the integration with respect
to $x$ gives
\begin{eqnarray}
 \int_0^\infty x g_{GIG(\delta,\gamma,\kappa)} dx &=& 
 \max(0, \kappa) \int_0^\infty e^{-x\frac{\gamma^2}{2}} +
 \int_0^\infty \int_0^\infty \frac{e^{-y(x + \gamma^2/2)}}{\pi^2 y  [J_{|\kappa|}^2 ( \sqrt{x}) + Y_{|\kappa|}^2(\sqrt{x})
  ] } dy dx  \\
&=& \frac{2}{\gamma^2} \max(0,\kappa) + \int_0^\infty \frac{1}{\pi (y
  + \gamma^2/2)y  [J_{|\kappa|}^2 ( \sqrt{x}) + Y_{|\kappa|}^2(\sqrt{x})
  ] } dy.
\end{eqnarray}
We now use the change of variables $y \rightarrow 2 \delta^2 y$ in the integral above and use the function
$g_\nu$ defined in (\ref{e:defofg}) to obtain
\begin{equation}
 \int_0^\infty \frac{1}{\pi (y
  + \gamma^2/2)y  [J_{|\kappa|^2} ( \sqrt{x}) + Y_{|\kappa|}^2(\sqrt{x})
  ] } dy = \delta^2 \int_0^\infty \frac{1}{y + \delta^2 \gamma^2}
g_{|\kappa|}(y) dy.
\end{equation}
Now, we apply the integral representation given in \cite{eberlein:2002}, formula (5.2) 
 to obtain
\begin{equation}
\delta^2 \int_0^\infty \frac{1}{y + \delta^2 \gamma^2}
g_{|\kappa|}(y) dy = \frac{\delta K_{|\kappa| - 1}(\delta
  \gamma)}{\gamma K_{|\kappa|}(\delta \gamma)}.
\end{equation}
Thus, (\ref{e:showthis1}) follows if we can now show
\begin{equation}\label{e:almost}
\frac{\delta K_{1+\kappa}( \gamma \delta) }{\gamma K_{\kappa}(\gamma
    \delta)} =  \frac{2}{\gamma^2} \max(0,\kappa) +  \frac{\delta K_{|\kappa| - 1}(\delta
  \gamma)}{\gamma K_{|\kappa|}(\delta \gamma)}, \quad \kappa \in
\mathbb{R}, \delta,\gamma > 0.
\end{equation}
For this, we require the following two properties of the Bessel function
$K_{\nu}$, (see for instance, \cite{gray:1966}, formulas (3.15) and (3.22))
\begin{eqnarray}
K_{\nu}(x) &=& K_{-\nu}(x), \quad \quad x \geq 0, \nu \in \mathbb{R} \label{e:symmetric} \\
x K_{\nu+2}(x)  &=&  x K_{\nu}(x) + 2 (1 + \nu) K_{1+\nu}(x) \quad x
\geq 0, \nu \in \mathbb{R}. \label{e:recursive}
\end{eqnarray}

For $\kappa \leq 0$, (\ref{e:almost}) follows immediately from
(\ref{e:symmetric}).  For $\kappa > 0$, (\ref{e:recursive}) gives
\begin{eqnarray}
\frac{2}{\gamma^2} \max(0,\kappa) +  \frac{\delta K_{|\kappa| - 1}(\delta
  \gamma)}{\gamma K_{|\kappa|}(\delta \gamma)} &=& \frac{2}{\gamma^2} \kappa +  \frac{\delta K_{\kappa - 1}(\delta
  \gamma)}{\gamma K_{\kappa}(\delta \gamma)}  \\
&=& \frac{2 \kappa K_{\kappa}(\delta \gamma) + \delta \gamma K_{\kappa
    -1}(\delta \gamma)}{\gamma^2 K_\kappa(\delta \gamma)} \\
&=& \frac{\delta K_{\kappa + 1}(\delta \gamma)}{\gamma
  K_{\kappa}(\delta \gamma)}.
\end{eqnarray}
This verifies (\ref{e:showthis1}), and hence finishes case 1.

\medskip

\noindent \textbf{Case 2:}  $\delta = 0, \gamma, \kappa > 0$.

This case is immediate, indeed,
\begin{equation}
i \lambda \frac{2 \kappa}{\gamma^2} + \int_0^\infty (e^{i \lambda x} -
1 -i \lambda x) \left(  \frac{\lambda}{x} e^{-\frac{\gamma^2}{2} x}
\right) dx = \int_0^\infty (e^{i \lambda x} - 1)
g_{GIG(0,\gamma,\kappa)}(x) dx.
\end{equation}

\medskip

\noindent \textbf{Case 3:}  $\gamma = 0, \delta>0, \kappa < 0$.

For this, we need to check that 
\begin{equation}
\delta^2 \int_0^\infty \frac{1 - e^{-x}}{x} g_{|\kappa|}(2 \delta^2 x)
  dx = \int_0^1 x g_{GIG(\delta,0,\kappa)}(x)dx
\end{equation}
This follows by changing the order of integration 
and using the definition of $g_\nu$:
\begin{eqnarray}
\int_0^1 x g_{GIG(\delta,0,\kappa)}(x)dx &=& \int_0^1 \int_0^\infty
\frac{e^{-y x }}{\pi y  [J_{|\kappa|}^2 ( \delta \sqrt{2y}) +
  Y_{|\kappa|}^2(\delta \sqrt{2y})
  ] } dy dx \\
&=&\int_0^\infty \frac{1 - e^{-y}}{y} \frac{1}{ \pi y  [J_{|\kappa|}^2
  ( \delta \sqrt{2 y}) + Y_{|\kappa|}^2(\delta\sqrt{2 y})
  ] } dy \\
&=& \delta^2 \int_0^\infty \frac{1 - e^{-y}}{y} g_{|\kappa|}(2 \delta^2 y)
  dy.
 \end{eqnarray}
This finishes the proof.

\end{pf}

}

\noindent Mark Veillette \& Murad Taqqu \\
\noindent \small Dept. of Mathematics \\
\noindent \small Boston University \\
\noindent \small 111 Cummington St. \\
\noindent \small Boston, MA 02215

\end{document}